\documentclass[a4paper,12pt]{amsart}
\usepackage[hiresbb]{graphicx}
\usepackage{amsmath,amssymb,latexsym}
\usepackage{tikz}
\usepackage{epic,eepic,eclarith}
\def\ArrowHeadSize{2}

\def\r{0.2}
\if00%
\def\c{0.9}
\newcommand{\Node}[4]{
\draw[very thick] (#1,#2) circle (#4) node  {$#3$};
}
\newcommand{\NodeRec}[4]{
\Mul{#4}{0.65}\acomo\Div{#4}{1.2}\rr
\draw[very thick] (#1,#2) node {$#3$}
   +(-\rr,-\rr) rectangle ++(\rr,\rr);
}
\newcommand{\Edge}[6]{
\Mul\c\c\s\Sub{1}\s\s\Sqroot\s\s
\Mul{#5}\s\zrone\Add{#1}\zrone\zrone
\Mul{#5}\c\zrtwo\Sub{#2}\zrtwo\zrtwo
\Mul{#5}\s\wrone\Add{#3}\wrone\wrone
\Mul{#5}\c\wrtwo\Add{#4}\wrtwo\wrtwo
\Add\zrone\wrone\zwrone\Div\zwrone{2}\zwrone
\Add\zrtwo\wrtwo\zwrtwo\Div\zwrtwo{2}\zwrtwo
\Add\zwrone{0}\zwronep
\Add\zwrtwo\acomo\zwrtwop
\Mul{#5}\s\zlone\Sub{#1}\zlone\zlone
\Mul{#5}\c\zltwo\Sub{#2}\zltwo\zltwo
\Mul{#5}\s\wlone\Sub{#3}\wlone\wlone
\Mul{#5}\c\wltwo\Add{#4}\wltwo\wltwo
\Add\zlone\wlone\zwlone\Div\zwlone{2}\zwlone
\Add\zltwo\wltwo\zwltwo\Div\zwltwo{2}\zwltwo
\Add\zwlone{0}\zwlonep
\Sub\zwltwo\acomo\zwltwop
\draw[very thick,#6] (\zrone,\zrtwo) -- (\wrone,\wrtwo);
\draw[very thick] (\zlone,\zltwo) -- (\wlone,\wltwo);
\draw[very thick,->,#6] (\zwrone,\zwrtwo) -- (\zwronep,\zwrtwop);
\draw[very thick,->] (\zwlone,\zwltwo) -- (\zwlonep,\zwltwop);
}
\newcommand{\EdgeOneUU}{
\def\Ac{.2}
\Sub{2}\r\Up
\Sub{1}\Ac\ArrowOne
\Add{1}\Ac\ArrowTwo
\draw[very thick,->] (0,\r) -- (0,1);
\draw[very thick,->] (0,\ArrowOne) -- (0,\ArrowTwo);
\draw[very thick] (0,\ArrowOne) -- (0,\Up);
}
\newcommand{\EdgeRec}[5]{
\Mul\c\c\s\Sub{1}\s\s\Sqroot\s\s
\Mul{#5}\s\zrone\Add{#1}\zrone\zrone
\Sub{#2}\rr\zrtwo%
\Mul{#5}\s\wrone\Add{#3}\wrone\wrone
\Mul{#5}\c\wrtwo\Add{#4}\wrtwo\wrtwo
\Add\zrone\wrone\zwrone\Div\zwrone{2}\zwrone
\Add\zrtwo\wrtwo\zwrtwo\Div\zwrtwo{2}\zwrtwo
\Add\zwrone{0}\zwronep
\Add\zwrtwo\acomo\zwrtwop
\Mul{#5}\s\zlone\Sub{#1}\zlone\zlone
\Mul{#5}\c\zltwo\Sub{#2}\zltwo\zltwo
\Mul{#5}\s\wlone\Sub{#3}\wlone\wlone
\Mul{#5}\c\wltwo\Add{#4}\wltwo\wltwo
\Add\zlone\wlone\zwlone\Div\zwlone{2}\zwlone
\Add\zltwo\wltwo\zwltwo\Div\zwltwo{2}\zwltwo
\Add\zwlone{0}\zwlonep
\Sub\zwltwo\acomo\zwltwop
\draw[very thick] (\zrone,\zrtwo) -- (\wrone,\wrtwo);
\draw[very thick] (\zlone,\zltwo) -- (\wlone,\wltwo);
\draw[very thick,->] (\zwrone,\zwrtwo) -- (\zwronep,\zwrtwop);
\draw[very thick,->] (\zwlone,\zwltwo) -- (\zwlonep,\zwltwop);
}
\newcommand{\EdgeCRec}[5]{
\def\r{#5}
\Mul\c\c\s\Sub{1}\s\s\Sqroot\s\s
\Mul\r\c\zuone\Add{#1}\zuone\zuone\Mul\r\s\zutwo\Add{#2}\zutwo\zutwo
\Mul\r\c\zdone\Add{#1}\zdone\zdone\Mul\r\s\zdtwo\Sub{#2}\zdtwo\zdtwo
\Sub{#3}\rr\zzuone
\Mul\r\s\rs
\Add{#4}\rs\zzutwo
\Sub{#4}\rs\zzdtwo
\Add\zuone\zzuone\zzzuone\Div\zzzuone{2}\zzzuone
\Add\zzzuone\acomo\zzzuoneacomo
\Sub\zzzuone\acomo\zzzuoneacomoD
\Mul\r\s\zlone\Sub{#1}\zlone\zlone
\Mul\r\c\zltwo\Sub{#2}\zltwo\zltwo
\Mul\r\s\wlone\Sub{#3}\wlone\wlone
\Mul\r\c\wltwo\Add{#4}\wltwo\wltwo
\Add\zlone\wlone\zwlone\Div\zwlone{2}\zwlone
\Add\zltwo\wltwo\zwltwo\Div\zwltwo{2}\zwltwo
\Add\zwlone{0}\zwlonep
\Sub\zwltwo\acomo\zwltwop
\draw[very thick] (\zuone,\zutwo) -- (\zzuone,\zzutwo);
\draw[very thick] (\zuone,\zdtwo) -- (\zzuone,\zzdtwo);
\draw[very thick,->] (\zzzuone,\zutwo) -- (\zzzuoneacomo,\zzutwo);
\draw[very thick,->] (\zzzuone,\zdtwo) -- (\zzzuoneacomoD,\zdtwo);
}
\newcommand{\EdgeCRecRLL}[5]{
\def\r{#5}
\Mul\c\c\s\Sub{1}\s\s\Sqroot\s\s
\Mul\r\c\zuone\Add{#1}\zuone\zuone\Mul\r\s\zutwo\Add{#2}\zutwo\zutwo
\Mul\r\c\zdone\Add{#1}\zdone\zdone\Mul\r\s\zdtwo\Sub{#2}\zdtwo\zdtwo
\Sub{#3}\rr\zzuone
\Mul\r\s\rs
\Add{#4}\rs\zzutwo
\Sub{#4}\rs\zzdtwo
\Add\zuone\zzuone\zzzuone\Div\zzzuone{2}\zzzuone
\Add\zzzuone\acomo\zzzuoneacomo
\Sub\zzzuone\acomo\zzzuoneacomoD
\Mul\r\s\zlone\Sub{#1}\zlone\zlone
\Mul\r\c\zltwo\Sub{#2}\zltwo\zltwo
\Mul\r\s\wlone\Sub{#3}\wlone\wlone
\Mul\r\c\wltwo\Add{#4}\wltwo\wltwo
\Add\zlone\wlone\zwlone\Div\zwlone{2}\zwlone
\Add\zltwo\wltwo\zwltwo\Div\zwltwo{2}\zwltwo
\Add\zwlone{0}\zwlonep
\Sub\zwltwo\acomo\zwltwop
\draw[very thick] (\zuone,\zutwo) -- (\zzuone,\zzutwo);
\draw[very thick] (\zuone,\zdtwo) -- (\zzuone,\zzdtwo);
\draw[very thick,->] (\zzzuone,\zutwo) -- (\zzzuoneacomo,\zzutwo);
\draw[very thick,->] (\zzzuone,\zdtwo) -- (\zzzuoneacomoD,\zdtwo);
\Add\zzzuone{.2}\zzzuoneD
\Add\zzzuoneacomoD{.2}\zzzuoneacomoDD
\draw[very thick,->] (\zzzuoneD,\zdtwo) -- (\zzzuoneacomoDD,\zdtwo);
}
\newcommand{\Acc}[1]{
\draw[very thick] (1,3) node {$#1$};
}

\newcommand{\CtwoLoop}[4]{
\def\r{0.3}\Mul\r{0.65}\acomo\def\c{0.9}\Div\r{1.2}\rr
\def\hh{0.8}\def\ll{1.5}\def\he{0.8}\Add\ll{0.2}\lll
\Div\he{3}\hehe
\begin{tikzpicture}[scale=1,>=stealth]
\useasboundingbox (-2,0) rectangle (2,-2);
\Node{-1}{0}{#1}{#3}
\LoopU{-1}{.1}{1.2}{}
\NodeRec{1}{0}{#2}{#3}
\LoopU{1}{.1}{1.2}{}
\EdgeCRec{-1}{0}{1}{0}{#3}  
\draw[very thick] (0,\hh) node {#4}; %
\end{tikzpicture}
}
\newcommand{\CtwoZXRLL}[4]{
\def\r{0.25}\Mul\r{0.65}\acomo\def\c{0.9}\Div\r{1.2}\rr
\def\hh{0.8}\def\ll{1.5}\def\he{0.8}\Add\ll{0.2}\lll
\Div\he{3}\hehe
\begin{tikzpicture}[scale=1,>=stealth]
\useasboundingbox (-2,0) rectangle (2,-2);
\Node{-1}{0}{#1}{#3}
\NodeRec{1}{0}{#2}{#3}
\LoopUDthree{-1}{0}{1}{}
\LoopURec{1}{0}{1.2}{}
\EdgeCRecRLL{-1}{0}{1}{0}{#3}  
\draw[very thick] (0,1.5) node {#4}; %
\end{tikzpicture}
}
\newcommand{\CtwoK}[8]{
\def\r{0.3}\Mul\r{0.65}\acomo\def\c{0.9}\Div\r{1.2}\rr
\def\hh{0.8}\def\ll{1.3}\def\he{#7}\Add\ll{0.2}\lll
\Div\he{3}\hehe
\begin{tikzpicture}[scale=#8,>=stealth]
\useasboundingbox (#5) rectangle (#6);
\draw[very thick] (-\ll,\he).. controls (-\lll,\hehe) and (-\lll,-\hehe)
 ..(-\ll,-\he);
\draw[very thick] (\ll,\he).. controls (\lll,\hehe) and (\lll,-\hehe)
 ..(\ll,-\he);
\Node{-1}{0}{#1}{#3}
\NodeRec{1}{0}{#2}{#3}
\EdgeCRec{-1}{0}{1}{0}{#3}  
\draw[very thick] (0,\hh) node {#4}; %
\end{tikzpicture}
}

\newcommand{\CtwoNKs}[8]{
\def\r{0.3}\Mul\r{0.65}\acomo\def\c{0.9}\Div\r{1.2}\rr
\def\hh{0.8}\def\ll{1.3}\def\he{#7}\Add\ll{0.2}\lll
\Div\he{3}\hehe
\begin{tikzpicture}[scale=#8,>=stealth]
\useasboundingbox (#5) rectangle (#6);
\draw[very thick] (-\ll,\he).. controls (-\lll,\hehe) and (-\lll,-\hehe)
 ..(-\ll,-\he);
\draw[very thick] (\ll,\he).. controls (\lll,\hehe) and (\lll,-\hehe)
 ..(\ll,-\he);
\Node{-1}{0}{#1}{#3}
\NodeRec{1}{0}{#2}{#3}
\draw[very thick] (0,\hh) node {#4}; %
\end{tikzpicture}
}

\newcommand{\CtwoLoopKs}[8]{
\def\r{0.3}\Mul\r{0.65}\acomo\def\c{0.9}\Div\r{1.2}\rr
\def\hh{0.8}\def\ll{1.5}\def\he{#7}\Add\ll{0.2}\lll
\Div\he{3}\hehe
\begin{tikzpicture}[scale=#8,>=stealth]
\useasboundingbox (#5) rectangle (#6);
\draw[very thick] (-\ll,\he).. controls (-\lll,\hehe) and (-\lll,-\hehe)
 ..(-\ll,-\he);
\draw[very thick] (\ll,\he).. controls (\lll,\hehe) and (\lll,-\hehe)
 ..(\ll,-\he);
\Node{-1}{0}{#1}{#3}
\NodeRec{1}{0}{#2}{#3}
\EdgeCRec{-1}{0}{1}{0}{#3}  
\draw[very thick] (0,\hh) node {#4}; %
\Mul{#3}{2.2}\rad
\Sub\rad{#3}\radr
\Add\rad\radr\radrr
\draw[very thick,red] (-1,\rad) circle (\radr);
\draw[very thick,red] (1,\rad) circle (\radr);
\draw[very thick,->,red] (-0.9,\radrr)-- (-1.1,\radrr); 
\draw[very thick,<-,red] (0.9,\radrr)-- (1.1,\radrr); 
\end{tikzpicture}
}

\newcommand{\ConeKs}[8]{
\def\r{0.3}\Mul\r{0.65}\acomo\def\c{0.9}\Div\r{1.2}\rr
\def\hh{0.8}\def\ll{0.5}\def\he{#7}\Add\ll{0.2}\lll
\Add\he\rr\hee
\Div\he{3}\hehe
\Div\hee{2}\hehee
\begin{tikzpicture}[scale=#8,>=stealth]
\useasboundingbox (#5) rectangle (#6);
\draw[very thick] (-\ll,\hee).. controls (-\lll,\hehee) and (-\lll,-\hehe)
 ..(-\ll,-\he);
\draw[very thick] (\ll,\hee).. controls (\lll,\hehee) and (\lll,-\hehe)
 ..(\ll,-\he);
\NodeRec{0}{0}{#1}{#3}
\LoopURec{0}{-.05}{1.2}{}
\draw[very thick] (0,\hh) node {#4}; %
\end{tikzpicture}
}
\newcommand{\ConeMs}[8]{
\def\r{0.3}\Mul\r{0.65}\acomo\def\c{0.9}\Div\r{1.2}\rr
\def\hh{0.8}\def\ll{0.5}\def\he{#7}\Add\ll{0.2}\lll
\Add\he\rr\hee
\Div\he{3}\hehe
\Div\hee{2}\hehee
\begin{tikzpicture}[scale=#8,>=stealth]
\useasboundingbox (#5) rectangle   (#6);
\draw[very thick] (-\ll,\hee).. controls (-\lll,\hehee) and (-\lll,-\hehe)
 ..(-\ll,-\he);
\draw[very thick] (\ll,\hee).. controls (\lll,\hehee) and (\lll,-\hehe)
 ..(\ll,-\he);
\Node{0}{0}{#1}{#3}
\LoopURec{0}{-.05}{1.2}{}
\draw[very thick] (0,\hh) node {#4}; %
\end{tikzpicture}
}
\newcommand{\ConeKK}[5]{
\def\r{0.3}\Mul\r{0.65}\acomo\def\c{0.9}\Div\r{1.2}\rr
\def\hh{0.8}\def\ll{0.5}\def\he{#5}\Add\ll{0.2}\lll
\Add\he\rr\hee
\Div\he{3}\hehe
\Div\hee{2}\hehee
\draw[very thick] (-\ll,\hee).. controls (-\lll,\hehee) and (-\lll,-\hehe)
 ..(-\ll,-\he);
\draw[very thick] (\ll,\hee).. controls (\lll,\hehee) and (\lll,-\hehe)
 ..(\ll,-\he);
\NodeRec{0}{0}{#1}{#3}
\LoopURec{0}{-.05}{1.2}{}
\draw[very thick] (0,\hh) node {#4}; %
}
\newcommand{\EdgeCRecrll}[5]{
\def\r{#5}
\Mul\c\c\s\Sub{1}\s\s\Sqroot\s\s
\Mul\r\c\zuone\Add{#1}\zuone\zuone\Mul\r\s\zutwo\Add{#2}\zutwo\zutwo
\Mul\r\c\zdone\Add{#1}\zdone\zdone\Mul\r\s\zdtwo\Sub{#2}\zdtwo\zdtwo
\Sub{#3}\rr\zzuone
\Mul\r\s\rs
\Add{#4}\rs\zzutwo
\Sub{#4}\rs\zzdtwo
\Add\zuone\zzuone\zzzuone\Div\zzzuone{2}\zzzuone
\Add\zzzuone\acomo\zzzuoneacomo
\Sub\zzzuone\acomo\zzzuoneacomoD
\Mul\r\s\zlone\Sub{#1}\zlone\zlone
\Mul\r\c\zltwo\Sub{#2}\zltwo\zltwo
\Mul\r\s\wlone\Sub{#3}\wlone\wlone
\Mul\r\c\wltwo\Add{#4}\wltwo\wltwo
\Add\zlone\wlone\zwlone\Div\zwlone{2}\zwlone
\Add\zltwo\wltwo\zwltwo\Div\zwltwo{2}\zwltwo
\Add\zwlone{0}\zwlonep
\Sub\zwltwo\acomo\zwltwop
\draw[very thick] (\zuone,\zutwo) -- (\zzuone,\zzutwo);
\draw[very thick] (\zuone,\zdtwo) -- (\zzuone,\zzdtwo);
\draw[very thick,->] (\zzzuone,\zutwo) -- (\zzzuoneacomo,\zzutwo);
\Sub\zzzuone{0.1}\zzzuoneT\Sub\zzzuoneacomoD{0.1}\zzzuoneacomoDT
\Add\zzzuone{0.1}\zzzuoneTT\Add\zzzuoneacomoD{0.1}\zzzuoneacomoDTT
\draw[very thick,->] (\zzzuoneTT,\zdtwo) -- (\zzzuoneacomoDTT,\zdtwo);
\draw[very thick,->] (\zzzuoneT,\zdtwo) -- (\zzzuoneacomoDT,\zdtwo);
}
\newcommand{\EdgeCCrrl}[5]{
\def\r{#5}\def\c{0.9}
\def\zx{#1}\def\zy{#2}\def\wx{#3}\def\wy{#4}
\Mul\c\c\s\Sub{1}\s\s\Sqroot\s\s
\Div\r{1.4142}\zzx\Add{#1}\zzx\zzx  
\Div\r{1.4142}\zzy\Add{#2}\zzy\zzy 
\Sub\zzx\zx\zzzx\Sub\zzy\zy\zzzy
\Mul\c\zzzx\czzzx\Mul\s\zzzy\szzzy 
\Mul\s\zzzx\szzzx\Mul\c\zzzy\czzzy 
\Add\czzzx\szzzy\zrx\Add\zrx\zx\zrx
\Sub\czzzy\szzzx\zry\Add\zry\zy\zry 
\Sub\czzzx\szzzy\zlx\Add\zlx\zx\zlx
\Add\szzzx\czzzy\zly\Add\zly\zy\zly
\Div\r{1.4142}\wwx\Sub{#3}\wwx\wwx  
\Div\r{1.4142}\wwy\Sub{#4}\wwy\wwy 
\Sub\wwx\wx\wwwx\Sub\wwy\wy\wwwy
\Mul\c\wwwx\cwwwx\Mul\s\wwwy\swwwy 
\Mul\s\wwwx\swwwx\Mul\c\wwwy\cwwwy 
\Sub\cwwwx\swwwy\wrx\Add\wrx\wx\wrx
\Add\cwwwy\swwwx\wry\Add\wry\wy\wry 
\Add\cwwwx\swwwy\wlx\Add\wlx\wx\wlx
\Sub\cwwwy\swwwx\wly\Add\wly\wy\wly
\Sub\wlx\zlx\zwlx\Div\zwlx{2}\zwlx\Add\zlx\zwlx\zwlx
\Sub\wly\zly\zwly\Div\zwly{2}\zwly\Add\zly\zwly\zwly
\def\ac{0.15}\def\acc{-0.01}
\Sub\zwlx{\acc}\zwlxA\Sub\zwly{\acc}\zwlyA
\Add\zwlx{\ac}\zwlxB\Add\zwly{\ac}\zwlyB
\Sub\zwlxA{0.02}\zwlxAA\Sub\zwlyA{0.02}\zwlyAA
\Sub\zwlxB{0.02}\zwlxBB\Sub\zwlyB{0.02}\zwlyBB
\draw[very thick,->] (\zlx,\zly) -- (\zwlxA,\zwlyA);
\draw[very thick,->] (\zwlxAA,\zwlyAA) -- (\zwlxB,\zwlyB);
\draw[very thick] (\zwlxBB,\zwlyBB) -- (\wlx,\wly);
\Sub\wrx\zrx\zwrx\Div\zwrx{2}\zwrx\Add\zrx\zwrx\zwrx
\Sub\wry\zry\zwry\Div\zwry{2}\zwry\Add\zry\zwry\zwry
\def\ac{0.15}\def\acc{0.1}
\Sub\zwrx{\acc}\zwrxA\Sub\zwry{\acc}\zwryA
\Add\zwrx{\ac}\zwrxB\Add\zwry{\ac}\zwryB
\Sub\zwrxA{0.02}\zwrxAA\Sub\zwryA{0.02}\zwryAA
\Sub\zwrxB{0.02}\zwrxBB\Sub\zwryB{0.02}\zwryBB
\draw[very thick] (\zrx,\zry) -- (\zwrxA,\zwryA);
\draw[very thick,<-] (\zwrxA,\zwryA) -- (\wrx,\wry);
}
\newcommand{\EdgeCC}[5]{
\def\r{#5}\def\c{0.9}
\def\zx{#1}\def\zy{#2}\def\wx{#3}\def\wy{#4}
\Mul\c\c\s\Sub{1}\s\s\Sqroot\s\s
\Div\r{1.4142}\zzx\Sub{#1}\zzx\zzx  
\Div\r{1.4142}\zzy\Add{#2}\zzy\zzy 
\Sub\zzx\zx\zzzx\Add\zzy\zy\zzzy
\Div\r{1.4142}\wwx\Add{#3}\wwx\wwx  
\Div\r{1.4142}\wwy\Sub{#4}\wwy\wwy 
\Add\wwx\wx\wwwx\Sub\wwy\wy\wwwy
\Sub\wx\zx\zwx\Div\zwx{2}\zwx\Add\zx\zwx\zwx
\Sub\wy\zy\zwy\Div\zwy{2}\zwy\Add\zy\zwy\zwy
\def\acx{-0.15}\def\acy{0.15}
\def\accx{-0.01}\def\accy{0.01}
\Add\zwx\acx\zwxB\Add\zwy\acy\zwyB
\Sub\zwx\accx\zwxA\Sub\zwy\accy\zwyA
\draw[very thick] (\zzx,\zzy) -- (\zwx,\zwy);
\draw[very thick,->] (\zwxA,\zwyA) -- (\zwxB,\zwyB);
\draw[very thick] (\zwxA,\zwyA) -- (\wwx,\wwy);
}
\newcommand{\EdgeCRecCurve}[5]{
\def\r{#5}\def\c{0.9}
\def\zx{#1}\def\zy{#2}\def\wx{#3}\def\wy{#4}
\Mul\c\c\s\Sub{1}\s\s\Sqroot\s\s
\Div\r{1.4142}\zzx\Sub{#1}\zzx\zzx  
\Div\r{1.4142}\zzy\Sub{#2}\zzy\zzy 
\Div\r{1.4142}\wwx\Sub{#3}\wwx\wwx  
\Div\r{1.4142}\wwy\Sub{#4}\wwy\wwy 
\draw[very thick] (\zzx,\zzy) .. 
controls (-1,-0.7) and (0.3,-1.7) ..(\wwx,\wwy);
\draw[very thick,->] (.5,-.93)-- (.52,-.93);
}
\newcommand{\LoopU}[4]{
\def\r{0.2}\Mul\r{#3}\rr
\Add{#2}\r\Loopr\Add\rr\Loopr\Loopr
\draw[very thick,#4] (#1,\Loopr) circle (\rr);
\Add{#1}{.1}\ArrowFrom
\Sub{#1}{.15}\ArrowTo
\Add\Loopr\rr\ArrowY
\draw[very thick,->,#4] (\ArrowFrom,\ArrowY)-- (\ArrowTo,\ArrowY); 
}
\newcommand{\LoopURec}[4]{
\def\r{0.2}\Mul\r{#3}\rr
\Add{#2}\r\Loopr\Add\rr\Loopr\Loopr
\Add\Loopr{.03}\Looprr
\draw[very thick,#4] (#1,\Looprr) circle (\rr);
\Add{#1}{.1}\ArrowFrom
\Sub{#1}{.15}\ArrowTo
\Add\Looprr\rr\ArrowY
\draw[very thick,->,#4] (\ArrowFrom,\ArrowY)-- (\ArrowTo,\ArrowY); 
}
\newcommand{\LoopUDthree}[4]{
\def\cc{0.75}\Mul\cc\cc\ss\Sub{1}\ss\ss\Sqroot\ss\ss
\def\r{0.25}\Mul\r{#3}\rr
\Add{#2}\r\Loopr\Add\rr\Loopr\Loopr
\draw[very thick,#4] (#1,\Loopr) circle (\rr);
\Add{#1}{.06}\ArrowFrom
\Sub{#1}{.13}\ArrowTo
\Add\Loopr\rr\ArrowY
\draw[very thick,->,#4] (\ArrowFrom,\ArrowY)-- (\ArrowTo,\ArrowY); 
\Mul\rr\ss\Arrx\Add{#1}\Arrx\Arrx
\Mul\rr\cc\Arry\Add\Loopr\Arry\Arry
\Mul\rr\ss\Arlx\Sub{#1}\Arlx\Arlx
\Mul\rr\cc\Arly\Add\Loopr\Arly\Arly
\draw[very thick,->] (\Arrx,\Arry) -- ++(-.1,.1); 
\draw[very thick,->] (\Arlx,\Arly) -- ++(-.1,-.1); 
}
\newcommand{\LoopUDtwo}[4]{
\def\cc{0.81}\Mul\cc\cc\ss\Sub{1}\ss\ss\Sqroot\ss\ss
\def\r{0.25}\Mul\r{#3}\rr
\Add{#2}\r\Loopr\Add\rr\Loopr\Loopr
\draw[very thick,#4] (#1,\Loopr) circle (\rr);
\Add{#1}{.06}\ArrowFrom
\Sub{#1}{.13}\ArrowTo
\Add\Loopr\rr\ArrowY
\Mul\rr\ss\Arrx\Add{#1}\Arrx\Arrx
\Mul\rr\cc\Arry\Add\Loopr\Arry\Arry
\Mul\rr\ss\Arlx\Sub{#1}\Arlx\Arlx
\Mul\rr\cc\Arly\Add\Loopr\Arly\Arly
\draw[very thick,->] (\Arrx,\Arry) -- ++(-.1,.1); 
\draw[very thick,->] (\Arlx,\Arly) -- ++(-.1,-.1); 
}
\newcommand{\LoopL}[4]{
\def\r{0.2}\Mul\r{#3}\rr
\Sub{#1}\r\Loopr\Sub\Loopr\rr\Loopr
\draw[very thick,#4] (\Loopr,#2) circle (\rr);
\Add{#2}{.1}\ArrowFrom
\Sub{#2}{.15}\ArrowTo
\Sub\Loopr\rr\ArrowX
\draw[very thick,->,#4] (\ArrowX,\ArrowFrom)-- (\ArrowX,\ArrowTo); 
}
\newcommand{\LoopD}[4]{
\def\r{0.2}\Mul\r{#3}\rr
\Sub{#2}\r\Loopr\Sub\Loopr\rr\Loopr
\draw[very thick,#4] (#1,\Loopr) circle (\rr);
\Sub{#1}{.1}\ArrowFrom
\Add{#1}{.15}\ArrowTo
\Sub\Loopr\rr\ArrowY
\draw[very thick,->,#4] (\ArrowFrom,\ArrowY)-- (\ArrowTo,\ArrowY); 
}

\newcommand{\Dfive}{
\def\R{0.2}
\Node{-1}{1}{3}{\R}
\Node{0}{2}{1}{\R}
\Node{0}{0}{2}{\R}
\NodeRec{2}{2}{1}{\R}
\NodeRec{2}{0}{2}{\R}
\Edge{0}{2}{0}{0}{\R}{}
\EdgeRec{2}{2}{2}{0}{\R}
\EdgeCRecrll{0}{2}{2}{2}{\R}
\EdgeCRecrll{0}{0}{2}{0}{\R}
\EdgeCCrrl{-1}{1}{0}{2}{\R}
\EdgeCC{0}{0}{-1}{1}{\R}
\EdgeCRecCurve{-1}{1}{2}{0}{\R}
\LoopU{0}{2}{1.2}{}
\LoopL{-1}{1}{1.2}{}
\LoopD{0}{0}{1.2}{}
}
\newcommand{\Dfivethree}{
\def\R{0.2}
\Node{0}{2}{1}{\R}
\Node{0}{0}{2}{\R}
\NodeRec{2}{2}{1}{\R}
\Edge{0}{2}{0}{0}{\R}{}
\EdgeCRecrll{0}{2}{2}{2}{\R}
\LoopU{0}{2}{1.2}{}
\LoopD{0}{0}{1.2}{}
}
\newcommand{\DfiveT}{
\begin{tikzpicture}[scale=1.2,>=stealth]
\useasboundingbox (-1,-1) rectangle (1,1);
\Dfive\Acc{X}
\end{tikzpicture}
}
\newcommand{\KakkoLL}[5]{
\begin{tikzpicture}[scale=#5,>=stealth] 
\Mul{#2}{0.1}\ll
\def\he{#2}\Mul\ll{3}\lll\Div\he{2}\hehe
\useasboundingbox (#3) rectangle (#4);
\draw[very thick] (-\ll,\he).. controls (-\lll,\hehe) and (-\lll,-\hehe)
 ..(-\ll,-\he);
\end{tikzpicture}
\hspace{-11mm}
}

\newcommand{\KakkoRs}[5]{
\begin{tikzpicture}[scale=#5,>=stealth] 
\Mul{#2}{0.1}\ll
\def\he{#2}\Mul\ll{3}\lll\Div\he{2}\hehe
\useasboundingbox (#3) rectangle (#4);
\draw[very thick] (\ll,\he).. controls (\lll,\hehe) and (\lll,-\hehe)
 ..(\ll,-\he);
\end{tikzpicture}
\hspace{-11mm}
}
\newcommand{\LKakko}[2]{
\def\Kakkol{-#2}\def\KakkouD{#1}
\Add{2}\KakkouD\Kakkou\Sub{0}\KakkouD\Kakkod 
\def\KD{0.4}
\Sub\Kakkol\KD\Kakkoll\Sub\Kakkou\KD\Kakkouu
\Add\Kakkod\KD\Kakkodd
\draw[very thick] 
(\Kakkol,\Kakkou).. controls (\Kakkoll,\Kakkouu) 
and  (\Kakkoll,\Kakkodd) .. (\Kakkol,\Kakkod);
\hspace*{0mm}
}
\newcommand{\RKakko}[3]{
\def\Kakkol{#2}\def\KakkouD{#1}
\Add{2}\KakkouD\Kakkou
\Sub{0}\KakkouD\Kakkod 
\def\KD{0.4}
\Add\Kakkol\KD\Kakkoll
\Sub\Kakkou\KD\Kakkouu
\Add\Kakkod\KD\Kakkodd
\hspace*{#3mm}
\draw[very thick] 
(\Kakkol,\Kakkou).. controls (\Kakkoll,\Kakkouu) 
and  (\Kakkoll,\Kakkodd) .. (\Kakkol,\Kakkod);
}
\if00
\def\c{0.9}
\newcommand{\Cfour}{
\def\R{0.2}
\Node{0}{2}{1}{\R}
\Node{0}{0}{2}{\R}
\NodeRec{2}{2}{1}{\R}
\NodeRec{2}{0}{2}{\R}
\Edge{0}{2}{0}{0}{\R}{}
\EdgeRec{2}{2}{2}{0}{\R}
\EdgeCRec{0}{2}{2}{2}{\R}
\EdgeCRec{0}{0}{2}{0}{\R}
}
\def\marusikaku{-5pt}
\def\maruk#1#2{
\hspace{-1mm}\raisebox{\marusikaku}{\def\r{0.2} \tikz \draw[thick,#2] (0,0) circle (\r) node {#1} ;}
} 
\def\sikakuk#1#2{
\hspace{-.5mm}\raisebox{\marusikaku}{\def\r{0.2}\tikz \draw[thick,#2] (0,0) node {#1}
   +(-\r,-\r) rectangle ++(\r,\r);} 
}
\newcommand{\marukk}[3]{
\Mul{#3}{12}\zurasu\Mul{#3}{20}\zurasuu
\hspace{-\zurasu mm}\raisebox{\marusikaku}{
\tikz\draw[thick,#2] (0,0) circle (#3) node {#1} ;
\hspace{-\zurasuu mm}
}
}
\newcommand{\marukkrotate}[3]{
\Mul{#3}{12}\zurasu\Mul{#3}{20}\zurasuu
\hspace{-\zurasu mm}\raisebox{\marusikaku}{
\tikz\shadedraw[inner color=white,outer color=yellow,draw=blue]
(0,0)  node {#1} [rotate=45] +(-.2,-.2) rectangle ++(.2,.2);
\hspace{-\zurasuu mm}
}
}
\def\sikakukk#1#2#3{\Sub{#3}{0.02}\sik
\Mul{#3}{12}\zurasu\Mul{#3}{20}\zurasuu
\hspace{-\zurasu mm}\raisebox{\marusikaku}{
\begin{tikzpicture}[scale=1,>=stealth]
 \draw[thick,#2] (0,0) node {#1}
   +(-\sik,-\sik) rectangle ++(\sik,\sik);
\end{tikzpicture}
\hspace{-\zurasuu mm}
} 
}
\def\Petersentwo#1#2#3#4#5#6{
\def\zx{#1}\def\zy{#2}\def\r{#3}\def\x{#4}
\def\ax{\zx}\Mul\x{1}\ay\Add\zy\ay\ay
\Mul\x{-0.95}\bx\Add\zx\bx\bx\Mul\x{0.31}\by\Add\zy\by\by
\Mul\x{-0.59}\cx\Add\zx\cx\cx\Mul\x{-0.81}\cy\Add\zy\cy\cy
\Mul\x{0.59}\dx\Add\zx\dx\dx\Mul\x{-0.81}\dy\Add\zy\dy\dy
\Mul\x{0.95}\ex\Add\zx\ex\ex\Mul\x{0.31}\ey\Add\zy\ey\ey
\begin{tikzpicture}[scale=1,>=stealth]
\draw[very thick,#6] (\ax,\ay) circle (\r) node (a) {$2$};
\draw[very thick] (\bx,\by) node (b) {$1$};
\draw[very thick] (\bx,\by)   +(-\r,-\r) rectangle ++(\r,\r);
\draw[very thick] (\cx,\cy)  node (c) {$3$};
\draw[very thick] (\cx,\cy)   +(-\r,-\r) rectangle ++(\r,\r);
\draw[very thick] (\dx,\dy)  node (d) {$4$};
\draw[very thick] (\dx,\dy)   +(-\r,-\r) rectangle ++(\r,\r);
\draw[very thick] (\ex,\ey)  node (e) {$5$};
\draw[very thick] (\ex,\ey)   +(-\r,-\r) rectangle ++(\r,\r);
 \draw[very thick] (a) -- (c) -- (e) -- (b) -- (d) -- (a);
\Mul\x{2}\x
\def\ax{\zx}\Mul\x{1}\ay\Add\zy\ay\ay
\Mul\x{-0.95}\bx\Add\zx\bx\bx\Mul\x{0.31}\by\Add\zy\by\by
\Mul\x{-0.59}\cx\Add\zx\cx\cx\Mul\x{-0.81}\cy\Add\zy\cy\cy
\Mul\x{0.59}\dx\Add\zx\dx\dx\Mul\x{-0.81}\dy\Add\zy\dy\dy
\Mul\x{0.95}\ex\Add\zx\ex\ex\Mul\x{0.31}\ey\Add\zy\ey\ey
\shadedraw[inner color=white,outer color=yellow,draw=blue] 
[rotate=45] (1.18,1.17)  rectangle +(.42,.42);
\draw[very thick,#5] (\ax,\ay)   node (aa) {$1$};
\draw[very thick,#6] (\bx,\by) circle (\r)  node (bb) {$1$};
\draw[very thick] (\cx,\cy)   +(-\r,-\r) rectangle ++(\r,\r);
\draw[very thick] (\cx,\cy) node (cc) {$2$};
\draw[very thick] (\dx,\dy)   +(-\r,-\r) rectangle ++(\r,\r);
\draw[very thick] (\dx,\dy) node (dd) {$6$};
\draw[very thick,#6] (\ex,\ey)   circle (\r)  node (ee) {$3$};
 \draw[very thick] (a) -- (aa); 
 \draw[very thick] (b) -- (bb); 
 \draw[very thick] (c) -- (cc); 
 \draw[very thick] (d) -- (dd); 
 \draw[very thick] (e) -- (ee);
 \draw[very thick] (aa) -- (bb) -- (cc) -- (dd) -- (ee) -- (aa);
\end{tikzpicture}
}
\def\Petersenone#1#2#3#4#5{
\def\zx{#1}\def\zy{#2}\def\r{#3}\def\x{#4}
\def\ax{\zx}\Mul\x{1}\ay\Add\zy\ay\ay
\Mul\x{-0.95}\bx\Add\zx\bx\bx\Mul\x{0.31}\by\Add\zy\by\by
\Mul\x{-0.59}\cx\Add\zx\cx\cx\Mul\x{-0.81}\cy\Add\zy\cy\cy
\Mul\x{0.59}\dx\Add\zx\dx\dx\Mul\x{-0.81}\dy\Add\zy\dy\dy
\Mul\x{0.95}\ex\Add\zx\ex\ex\Mul\x{0.31}\ey\Add\zy\ey\ey
\begin{tikzpicture}[scale=1,>=stealth]
\draw[very thick,#5] (\ax,\ay) circle (\r) node (a) {$1$};
\draw[very thick,#5] (\bx,\by) circle (\r) node (b) {$2$};
\draw[very thick,#5] (\cx,\cy) circle (\r) node (c) {$3$};
\draw[very thick,#5] (\dx,\dy) circle (\r) node (d) {$4$};
\draw[very thick,#5] (\ex,\ey) circle (\r) node (e) {$5$};
 \draw[very thick] (a) -- (c) -- (e) -- (b) -- (d) -- (a);
\Mul\x{2}\x
\def\ax{\zx}\Mul\x{1}\ay\Add\zy\ay\ay
\Mul\x{-0.95}\bx\Add\zx\bx\bx\Mul\x{0.31}\by\Add\zy\by\by
\Mul\x{-0.59}\cx\Add\zx\cx\cx\Mul\x{-0.81}\cy\Add\zy\cy\cy
\Mul\x{0.59}\dx\Add\zx\dx\dx\Mul\x{-0.81}\dy\Add\zy\dy\dy
\Mul\x{0.95}\ex\Add\zx\ex\ex\Mul\x{0.31}\ey\Add\zy\ey\ey
\draw[very thick] (\ax,\ay)   +(-\r,-\r) rectangle ++(\r,\r);
\draw[very thick] (\ax,\ay) node (aa) {$1$};
\draw[very thick] (\bx,\by)   +(-\r,-\r) rectangle ++(\r,\r);
\draw[very thick] (\bx,\by) node (bb) {$2$};
\draw[very thick] (\cx,\cy)   +(-\r,-\r) rectangle ++(\r,\r);
\draw[very thick] (\cx,\cy) node (cc) {$3$};
\draw[very thick] (\dx,\dy)   +(-\r,-\r) rectangle ++(\r,\r);
\draw[very thick] (\dx,\dy) node (dd) {$4$};
\draw[very thick] (\ex,\ey)   +(-\r,-\r) rectangle ++(\r,\r);
\draw[very thick] (\ex,\ey) node (ee) {$5$};
 \draw[very thick] (a) -- (aa); 
 \draw[very thick] (b) -- (bb); 
 \draw[very thick] (c) -- (cc); 
 \draw[very thick] (d) -- (dd); 
 \draw[very thick] (e) -- (ee);
 \draw[very thick] (aa) -- (bb) -- (cc) -- (dd) -- (ee) -- (aa);
\end{tikzpicture}
}

\newcommand{\vv}{\overline{v}}
\newcommand{\References}[7]{[#1]\ #2,\  #3,\ #4 {\bf #5} (#7)\ #6.}

\newcommand{\zero}{\mathbf{0}}

\newcommand{\vt}[2]{v^{#1}_{#2}}
\newcommand{\charp}[1]{\phi\left(#1,\,x\right)}

\def\ddm{-}
\def\sumg{+}
\newcommand{\CC}{\mathbb{C}}

\newcommand{\ZZ}{\mathbb{Z}}
\newcommand{\EM}{\emph{Enseign. Math.}\ }

\newcommand{\yearb}[1]{}
\newcommand{\bbf}[1]{{#1},}%

\theoremstyle{plain}
\newtheorem{thm}{Theorem}[section]
\newtheorem{lem}[thm]{Lemma}
\newtheorem{cor}[thm]{Corollary}
\newtheorem{pro}[thm]{Proposition}

\theoremstyle{definition}
\newtheorem{defn}[thm]{Definition}
\newtheorem{rem}[thm]{Remark}
\newtheorem{ex}[thm]{Example}

\newcommand{\quotient}[1]{\hspace{-.5mm}/\hspace{-0mm}#1}
\newcommand{\Direct}[1]{\mathrm{SD}(#1)}
\def\sita#1{\kern-,1667em\raisebox{-.5ex}{\footnotesize $\,|#1$}\kern-.125em}
\def\sitaa{\kern-,1667em\raisebox{-.5ex}{$|$}\kern-.125em}
\newcommand{\vertex}[2]{v^{#1}_{#2}}
\newcommand{\vectortwo}[2]{\vertex{#1}{#2}}

\newcommand{\vi}{\{\overline{v}_i\}_{i=1}^r}
\newcommand{\smvi}[1]{#1\ssetminus \vi}
\newcommand{\ssetminus}{\hspace{-.8mm}\setminus\hspace{-1.2mm}}

\newcommand{\delG}[1]{#1\backslash\pi}
\newcommand{\DD}{\overline{D}}
\newcommand{\Pc}{P}

\newcommand{\dindex}[5]
{
#1_{
\mbox{\tiny $\vectortwo{#2}{#3}\vectortwo{#4}{#5}$}
   }
}  
\newcommand{\dindexonetwo}[4]
{
#1_{
\mbox{\tiny $#2\vectortwo{#3}{#4}$}
    }
}  
\newcommand{\dindextwoone}[4]
{
#1_{
\mbox{\tiny $\vectortwo{#2}{#3}#4$}
   }
}  
\newcommand{\one}{\bold{1}}
\newcommand{\tright}{\triangleright}

\newcommand{\Mat}[3]{\mathrm{Mat\,}(#1\times #2;\,#3)}

\newcommand{\dzi}{d^Z_i(X)}
\newcommand{\ddz}{d^Z}
\newcommand{\Dz}{D^Z(X)}

\newcommand{\CX}{A\left(\delG{X}\right)}
\newcommand{\rVd}[1]{#1\sita{V'\times V'}}

\begin{document}
{\Large
\begin{center}
Characteristic polynomials and zeta functions of equitably partitioned graphs
\end{center}
}
\bigskip
\vspace*{8mm}
\begin{center}
\large 
Osamu Kada\\
\end{center}
\begin{center}
\small
Part-time Lecturer, 
Faculty of Science and Engineering, \\
Hosei University, 
Koganei, Tokyo 184-8584, Japan
\end{center}
\begin{center}
\small
e-mail: \texttt{osamu.kada.77@hosei.ac.jp}
\end{center}
\vspace*{8mm}

\begin{center}
November, 2019
\end{center}
\vspace*{18mm}

\hspace{-4mm}{\bf Abstract}\quad 
 Let $\pi=\{V_1,\dots,V_r\}$ be an equitable partition of the vertex set of 
a directed graph (digraph) $X$. 
 It is well known that 
 the characteristic polynomial $\phi(X/\pi,x)$ of a quotient 
 graph $X/\pi$ divides that of $X$, but the remainder part 
is not well investigated.   
 In this paper, we define {\em a deletion graph}  
$\delG{X}$ over an equitable partition 
$\pi$, which is a signed directed graph  
defined for a fixed set of deleting vertices 
 $\{\overline{v}_i\in V_i, i=1,\cdots,r\}$,  
 and give a similarity transformation exchanging the adjacency matrix $A(X)$ 
 which is compatible with the equitable partition for 
 a block triangular matrix whose diagonal blocks are the adjacency matrix of 
the quotient graph and the deletion graph. 
 In fact, we show the result for more general matrices including 
 adjacency matrix of graphs, and 
 as  corollaries, we show the followings: 
 (i) a decomposition formula of the reciprocal 
of the Ihara-Bartholdi zeta function over an equitably partitioned 
 undirected graph into 
 the quotient graph part and the 
 deletion graph part, and  
 (ii) Chen and Chen's result ([CC17, Theorem 3.1])  on 
 the Ihara-Bartholdi zeta functions on generalized join graphs, 
 and (iii) Teranishi's result [Ter03, Theorem 3.3].  
  \bigskip\\
\if01
\def\r{2}
\def\rr{2.6}
 \begin{figure}[htbp]
\setlength{\unitlength}{.7mm}%
 	\begin{picture}(50,50)(0,0)%
\put(0,-20){
\put(20,80){\makebox(0,0){Figure 2.3}}
\put(30,-10){$2$'s best strategies}
\put(50,40){$1$'s best strategies}
\put(60,0){$x$}
\put(0,60){$y$}
\arrow\drawline(0,0)(0,50)
\arrow\drawline(0,0)(50,0)
\drawline(0,40)(40,40)(40,0)
\drawline(0,39.8)(20,39.8)
\drawline(0,40)(20,40)
\drawline(0,40.2)(20,40.2)
\drawline(0,40.4)(20,40.4)
\dottedline(20,40)(20,0)
\drawline(20,0)(40,0)
\drawline(20,-0.4)(40,-0.4)
\drawline(20,0.2)(40,0.2)
\drawline(20,0.4)(40,0.4)
\put(20,40){\blackcircle{\r}}
\put(19.6,0){\circle{\rr}}
\drawline(0,0)(0,18)
\drawline(0.2,0)(0.2,18)
\drawline(0.4,0)(0.4,18)
\put(0,19.5){\blackcircle{\r}}
\put(40,19.5){\circle{\rr}}
\drawline(40,20)(40,40)
\drawline(40.2,20)(40.2,40)
\drawline(40.4,20)(40.4,40)
\drawline(40.6,20)(40.6,40)
\dottedline(0,20)(40,20)
}
\end{picture}
 \end{figure}
\vspace*{20mm}
\Sqroot{2}\aaa
$\aaa$
\fi
\if00%
\section{Introduction}
 A partition $\pi=\{V_1,\dots,V_r\}$ of the vertex set $V(X)$ 
 of a directed graph  (digraph) $X$ (we allow multiple loops and edges) 
 is called {\em equitable} if 
 for each $i,j=1,\dots, r$, there is a integer $b_{ij}$ such that 
 each vertex in the cell $V_i$ has exactly $b_{ij}$ (oriented) edges  
 to vertices in the cell $V_j$. 
 The adjacency matrix of the (weighted) quotient graph over $\pi$ is defined by the matrix $(b_{ij})$ 
 and is denoted by $A(X\quotient{\pi})$.   
 It is well known that the characteristic polynomial 
 of the (weighted) quotient graph over $\pi$ (front divisor of $X$)
 divides the characteristic polynomial of the adjacency 
 matrix $A(X)$ of the graph ([GR01, 9.3]), but the remainder part 
is not well investigated.
 
 We define the {\em deletion graph} 
$\delG{X}$ over 
$\pi$  
 which is a signed directed graph  
defined for a fixed set of deleting 
 vertices $\{\overline{v}_i\in V_i, i=1,\cdots,r\}$,  
and give a similarity transformation exchanging the adjacency matrix $A(X)$ 
 which is compatible with the equitable partition for
 a block triangular matrix whose diagonal blocks are the adjacency matrices of 
the quotient graph and the deletion graph. 
 Hence, we have that 
 the remainder part 
 is the characteristic polynomial of 
 the adjacency matrix $A(\delG{X})$ of the deletion graph. 
 This answers the question posed by Deng and Wu [DW05, Section 5] of  
 whether we can associate any geometrical meaning 
 on the remainder part.

 Moreover, 
 we get the decomposition 
 formula for the characteristic  polynomials for the Laplacian matrix, 
 and the reciprocal of the Ihara-Bartholdi zeta function.  
 Since a generalized join (composition) (Schwenk [S74]) 
of regular graphs is a special case 
 of  an equitably partitioned graph, we get 
 the decomposition formula of the reciprocal of 
the Ihara-Bartholdi zeta function 
 of generalized joined graph by Chen and Chen [CC17] 
 as a corollary of our Theorem. 

 Having an equitably partitioned graph is equivalent 
 to having a covering projection 
([DSW07, Lemma 3.1]), 
and when  the graph is a covering of a voltage 
 assignment, this is equivalent to having a free action 
 (i.e. regular covering) ([GT77, Theorems 3,4], see also [DW05]).  
 There is much in the literature on the decomposition 
 of (a) the characteristic polynomial of the adjacency matrix of the graph  
 and (b) the reciprocal of the Ihara-Bartholdi zeta functions. 
 We list them.   
 On the topic (a), there are results 
 for graph covering with voltages in a finite group 
 by Mizuno and Sato [MS95, Theorem 1], [MS97, Theorem 1] 
(see also [KL92], [Sat99 Theorem 24], [KL01], [FKL04]); 
  for branched cover with  branch index 1  
by Deng and Wu [DW05, Theorem 4.2]  assuming a semi-free action on digraph;  
 for branched cover with 
branch index 1 by Deng, Sato  and Wu [DSW07, Theorem 6.4]. 
  On the topic (b), there are results for the reciprocal of 
 the (weighted) Ihara-zeta (Bartholdi-zeta) function of a 
 regular  ($g$-cyclic $\Gamma$-, or irregular) cover
 by Mizuno and Sato [MS01, Theorem 5], [MS02, Theorem 7],
 [MS04, Theorem 4], and by Sato [Sat06, Theorem 3], [Sat07, Theorem 4]. 

 When there is a symmetry (automorphism), so when the equitable partition 
 is the orbit partition, [BFW16], [FSSW17],[FSW18] 
 give a decomposition of any automorphism compatible matrix, which include 
 the adjacency matrix, the Laplacian matrix, etc.

 The remainder of the paper is organized as follows. 
 In section 2, we give basic facts on equitably partitioned 
 directed graphs. 
 In section 3, we define the deletion graph.
 In section 4, we  give our main theorem (Theorem \ref{zeta}) of 
 a similarity transformation exchanging the adjacency matrix $A(X)$ 
 for a block triangular matrix 
 whose diagonal blocks are the adjacency matrix of 
the quotient graph and the deletion graph, giving the decomposition 
 formula.   
 In Section 5 and 6 we give applications of the decomposition 
 formula to the reciprocal of the Ihara-Bartholdi zeta functions 
 of equitably partitioned graphs, 
 especially on generalized join graphs. 
\section{Equitable Partitions of Directed Graphs}
 For totally ordered  sets $U,V$ and a set $W$, we denote by $\Mat{U}{V}{W}$ 
  the set of matrices indexed by $U\times V$ whose components 
 are in $W$, that is, the set  $\{M:U\times V\to W\}$ of mappings 
 from $U\times V$ to $W.$
 For a square matrix  $M$,  $\phi(M,x)=\det(xI-M)$ 
 is the characteristic polynomial of $M$. 
 Let $X=(V(X),E(X))$ be a finite directed (multi)graph 
 with a set $V(X)$ of vertices and a set $E(X)$ of
 directed edges. 
 If $e\in E(X)$ implies  $\overline{e}\in E(X)$ 
 (here $\overline{e}$ is an inverse edge), 
 then $X$ can be considered as an undirected graph.   
 We allow $X$ to have multiple edges and multiple loops. 
 For $u,v\in V(X)$, we denote by $u\to v$ if there is an edge that 
 goes from $u$ to $v$, and $u\sim v$ if $u\to v$ and $v\to u$. 
 The adjacency matrix of $X$ is denoted by 
$A(X)$,  
 that is, $A(X)_{uv}$ is the number of edges from $u$ to $v$.
\begin{defn}
 Let $\pi=\{V_1,\dots,V_r\}$ be a partition of the vertex set $V(X)$ of 
 a directed graph $X$.
 For $i,j=1,\dots, r,$ for each vertex $u$ in the cell $V_i$,  
 if the number $b_{ij}$ of edges that goes from $u$ to the vertices in $V_j$ 
 does not depend on the choice of $u$,  
 we say that 
$\pi$ is an {\em equitable partition}.
 In this case, the multi-directed graph $X\quotient{\pi}$, called 
the {\em (weighted) quotient (or front divisor)}  of $X$ over $\pi$,  
 is such that the set of vertices is $\pi$, and there are $b_{ij}$ edges 
 from $V_i$ to $V_j$ ([GR01, 9.3]). 
The adjacency matrix of $X\quotient{\pi}$ is given by 
$A(X\quotient{\pi})=(b_{ij})\in \Mat{\pi}{\pi}{\ZZ_{+}}.$
\end{defn}
\begin{ex}\label{join}
\begin{enumerate}
 \item 
Let $X$ be the following $C_4$: 
 $\pi=\{V_1,V_2\}, V_i=\{v^i_1,v^i_2\},i=1,2,$ 
 $v^1_j=\marukk{j}{}{.25}\hspace{3mm},  
v^2_j=\sikakukk{j}{}{.25}\hspace{3mm}$. 
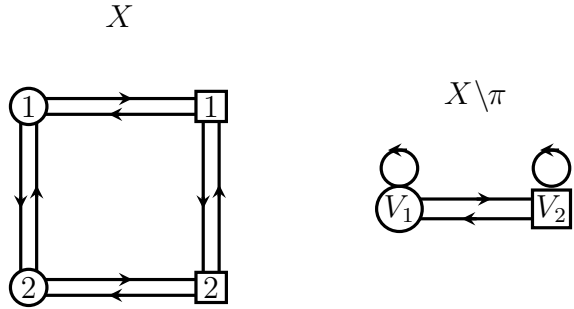
\vspace{15mm}

\begin{tikzpicture}[scale=1.2,>=stealth]
\useasboundingbox (-.8,-.8) rectangle (3,3);
\Cfour\Acc{X}
\end{tikzpicture}
\CtwoLoop{V_1}{V_2}{0.3}{\raisebox{15mm}{$X\ssetminus \pi$}}
\vspace{-15mm}

\begin{align*}
A(X)&=
 \begin{pmatrix}
 A_1& I_2 \\
I_2&A_1
\end{pmatrix},
 A_1= \begin{pmatrix}
 0& 1 \\
1&0
\end{pmatrix},
 A(X/\pi)=
\begin{pmatrix}
 1& 1\\
 1 & 1
\end{pmatrix},\\ 
\phi(A(X),x)&= x^2(x+2)(x-2),
\phi(A(X\quotient{\pi}),x)= x(x-2).
\end{align*}
\item
Next, consider the following directed graph.
 $\pi=\{V_1,V_2\}, V_1=\{v^1_j|\ j=1,2,3\}$, 
 $V_2=\{v^2_k|\ k=1,2\}$, 
 $v^1_j=\maruk{j}{}, v^2_k=\sikakuk{k}{}.$
\vspace{8mm}

\hspace*{11mm}
\DfiveT
\vspace*{.1mm}
\hspace*{21mm}
\CtwoZXRLL{V_1}{V_2}{0.25}{$X/\pi$}
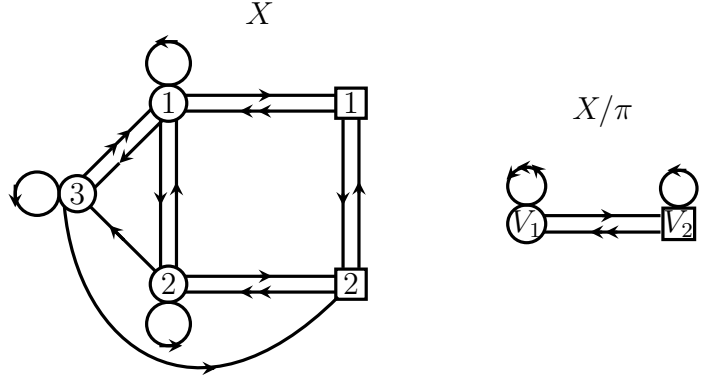
\vspace{-7mm}

\begin{align*}
  A(X)&=
\left(
\begin{array}{@{\,}ccc|cc@{\,}}
 1&1&1&1& 0 \\ 
1&1&1&0&1 \\
2&0&1 & 0&1 \\ \hline
2&0 & 0&0&1 \\ 
0&2&0&1& 0  
\end{array}
\right),\  
 A(X/\pi)=
\begin{pmatrix}
 3& 1\\
 2& 1
\end{pmatrix},\\ 
\phi(A(X),x)&=(x^2-4x+1)x(x-1)(x-2),\ 
 \phi(A(X\quotient{\pi}),x)=   
(x^2-4x+1). 
\end{align*}\vspace{-5mm}
\item
 We recall examples from [GR01; Section 9.3].
Let $X=J(5,2,0)$ be the Petersen graph.
(i)\ Let 
$\pi_1=\{V_1,V_2\}$, $V_i=\{\vectortwo{i}{k}|\ k=1,\dots,5\},\ i=1,2$,  
 $\vectortwo{1}{k}=$
\marukk{k}{red!96!white}{.2}, 
  $\vectortwo{2}{k}=$ \sikakuk{k}{black}.
\vspace{5mm}

\hspace{-8mm}
\Petersenone{0}{0}{0.2}{1}{red!96!white}
\vspace{-25mm}

$\hspace{48mm}
A(X)=
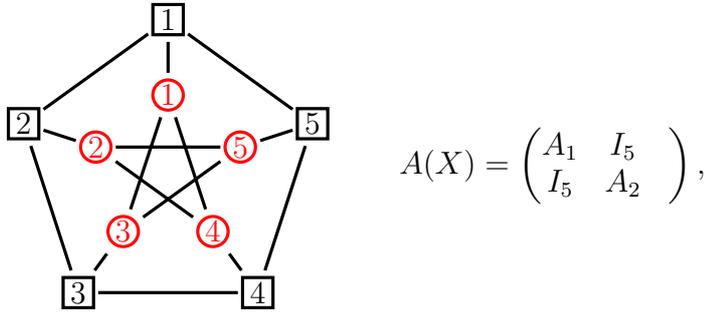
 \begin{pmatrix}
 A_1& I_5& \\
I_5&A_2
\end{pmatrix}, 
$
\vspace{13mm}

\begin{align*}
\hspace{-5mm}
A_1 & = \begin{pmatrix}
 0&0&1&1&0 \\
0&0&0&1&1\\
1&0&0&0&1\\
1&1&0&0&0\\
0&1&1&0&0
\end{pmatrix},
\hspace{2mm} 
A_2= \begin{pmatrix}
 0&1&0&0&1 \\
1&0&1&0&0\\
0&1&0&1&0\\
0&0&1&0&1\\
1&0&0&1&0
\end{pmatrix}, 
\hspace{2mm} A(X\quotient{\pi_1})=\begin{pmatrix}
 2 & 1 \\
 1 & 2
\end{pmatrix}, \\
&\phi(A(X),x)=(x-3)(x-1)^5(x+2)^4,\  \phi(A(X\quotient{\pi_1}),x)=(x-3)(x-1). 
\end{align*}
 (ii)\ Next consider the following distance partition $\pi_2=\{V_1,V_2,V_3\}$, 
 $V_1=\{\vectortwo{1}{1}\}$,  
 $V_2=\{\vectortwo{2}{k}|\ k=1,2,3\}$,  
 $V_3=\{\vectortwo{3}{l}|\ l=1,\dots,6\}$,   
  $\vertex{1}{1}=$\marukkrotate{1}{}{.2},
 $\vertex{2}{k}=\hspace{-1mm}$\marukk{k}{teal}{.2}, $k=1,2,3,$
  $\vertex{3}{l}=$\sikakuk{l}{black}, $l=1,\dots,6.$
 $V_i$ is the set of vetrices such that 
 their distances from $\vectortwo{1}{1}$ are the same $i-1$. 
\vspace{5mm}

\hspace{-5mm}
\Petersentwo{0}{0}{0.2}{1}{black}{teal}
\vspace{-50mm}

\begin{align*}
\hspace{65mm}
 A(X)&=
\left(
\begin{array}{@{\,}c|ccc|cccccc@{\,}}
 0&1&1&1& 0& 0& 0& 0& 0& 0 \\ \hline
1&0& 0& 0&1&1& 0& 0& 0& 0 \\
1&0& 0& 0& 0& 0&1&1& 0& 0 \\
1&0& 0& 0& 0& 0& 0& 0&1&1 \\ \hline
0&1&0& 0& 0& 0&0&1& 1& 0 \\
0&1&0& 0& 0& 0&1&0& 0& 1 \\
0&0&1& 0& 0& 1&0&0& 1& 0 \\
0&0&1& 0& 1& 0&0&0& 0& 1 \\
0&0& 0&1& 1& 0&0&1& 0& 0 \\
0&0& 0&1& 0& 1&0&1& 0& 0 
\end{array}
\right).
\end{align*}
$A(X\quotient{\pi_2})=\begin{pmatrix}
 0 & 3&0 \\
 1 & 0&2\\
0&1&2
\end{pmatrix}, \phi(A(X\quotient{\pi_2}),x)=(x-3)(x-1)(x+2).$ 
   
 For other examples of equitable partitions, see [GR01].
 \item  Let $X=H[X_1,\dots,X_r]$ be the generalized join (composition) 
 graph of 
$X_1,\dots,X_r$ determined by $H$ ([Sch74, 4.]), and assume $X_i$ are 
$k_i$-regular  
 ([Sch74, section 4], [CC17]). 
That is, $\# V(H)=r,V(X)=\cup_{i=1}^rV(X_i)$ and  
 for $\vertex{i}{k}\in V(X_i)$, $\vertex{j}{l}\in V(X_j)$, 
\begin{align}\label{Sch}
\dindex{A(X)}{i}{k}{j}{l}
=
\begin{cases}
 A(X_i)_{v^i_kv^i_l}& \text{ if }i=j,\\
A(H)_{ij}  & \text{ if }i\neq j.
\end{cases}  
\end{align}
 Then $\pi=\{V(X_1),\dots,V(X_r)\}$ is 
an equitable partition of $X$ by letting  
\begin{align}\label{quotient}
 A(X\quotient{\pi})_{ij}=  
\begin{cases}
   k_i & \text{if }i= j,\\
A(H)_{ij}n_j & \text{if }i\neq j, 
\end{cases}
\end{align}
 here $\ n_j=\# V_j$.  
\end{enumerate}  
\end{ex}
\if01
\begin{rem}
 Let $X_i$ be the subgraph of $X$ restricting 
to $V_i$. 
Let $A(X)\sita{V_i\times V_j}\in \Mat{V_i}{V_j}{\{0,1\}}$ be the restriction of $A(X)$ to $V_i\times V_j$, that is, the $V_i$-$V_j$ block matrix.
Then,  
$$
\pi \text{\ is equitable iff } \exists b_{ij}\in\mathbb{N} 
\text{ such that } 
 A(X)\sita{V_i\times V_j}\cdot\one=b_{ij}\cdot\one,
\ i,j=1,\dots,r, 
$$
 here $\one=
 \begin{pmatrix}1 \\ \vdots\\1\end{pmatrix}
\in \mathbb{Z}^{\# V_j}.$
 Motivating by this, we define as follows. 
\end{rem}
\fi
\begin{defn}
   Let $\pi=\{V_1,\dots,V_r\}$ be a partition of $V(X)$ of a digraph $X$
 and  let $M\in \Mat{V(X)}{V(X)}{\mathbb{C}}$.  
 We say that the pair $(M,\pi)$  is {\em equitable }
 if there exists $B=(b_{ij})\in \Mat{\pi}{\pi}{\mathbb{C}}$ 
 such that 
 $$M\sita{V_i\times V_j}\cdot\one =b_{ij}\cdot\one,\ i,j=1,\cdots,r,$$
 here $\one=
 \begin{pmatrix}1 \\ \vdots\\1\end{pmatrix}
\in \mathbb{Z}^{\# V_j}.$
 
  We denote by $M\quotient{\pi}:=B$, 
 and call it the {\em quotient matrix} of $M$ over the partition $\pi.$  
 Note that  $(A(X),\pi)$ is equitable  means that 
$\pi$ is an equitable partition, and  
$A(X)\quotient{\pi}=A(X\quotient{\pi})$.  
\end{defn}
\begin{rem}\label{equi}
 Let $\pi=\{V_1,\dots,V_r\}$ be an equitable partition of  $V(X)$, 
 and let 
  $M=\alpha A(X)+D$.  
 Here, $\alpha\in\mathbb{C}$ and $D$ is a diagonal matrix 
 such that 
$$D\sita{V_i\times V_i}=d_iI_i, 
\quad I_i\in \Mat{V_i}{V_i}{\{0,1\}} \text{ is the unit matrix.}
$$
 Then $(M,\pi)$ is equitable  with the quotient matrix  
\begin{equation}\label{equitable}
M\quotient{\pi}=
\alpha A(X\quotient{\pi})
+(\delta_{ij}d_i),
\end{equation}
 here $\delta_{ij}=1$ if $i=j$ and $=0$ if $i\neq j$. 
 Assume that $X$ is an undirected graph and 
 let   $D(X)$ be the degree matrix of $X$, that is, 
 the diagonal matrix such that 
  $D(X)_{uu}$ is the number of edges that go from $u$.
 Let $L(X)=-A(X)+D(X)$ (resp. $Q(X)=A(X)+D(X)$)  
be the Laplacian matrix (resp. the signless Laplacian  matrix)  
 of $X$, and let $Z_X(u,t)$ be the Ihara-Bartholdi zeta 
 function (see section 5). 
 Then by letting $(a) D=D(X)$, (i) $\alpha=-1,$ 
(ii) $\alpha=1,$ (b) $D=D^Z(X), \alpha=-t$,   
 we have  
\begin{align*}
 \text{(a) (i)\qquad  } L(X)\quotient{\pi}&=L(X\quotient{\pi})
=-A(X)\quotient{\pi})+(\delta_{ij}d_i(X)),\\
 \text{\quad (ii)\qquad  } Q(X)\quotient{\pi}&=Q(X\quotient{\pi})
=A(X\quotient{\pi})+(\delta_{ij}d_i(X)),\\
\text{(b)\hspace{5mm}  }Z_{X\quotient{\pi}}(u,t)^{-1}
&=s_1^{m-n}\det\left(-t A(X\quotient{\pi})+(\delta_{ij}d_i^Z(X))\right),
\end{align*}
 here $d_i(X)$ is the degree of vertices in $V_i$, 
 and $d_i^Z(X)=s_1+s_2d_i(X),$ ${s_1 =1-(1-u)^2t^2},$ $s_2=(1-u)t^2.$
\end{rem}
\section{The Deletion Graph over a partition of a graph}
\begin{defn}
 We define a signed directed graph $X=(V(X),\Direct{X})$ 
as follows: $\Direct{X}$
  is a set of a signed directed edges with its adjacency matrix 
 $A(X)\in\Mat{V(X)}{V(X)}{\mathbb{Z}}$ 
 whose components are integers. 
 That is, for $u,v\in V(X)$, if $A(X)_{uv}\ge 0$, 
 there are $A(X)_{uv}$ positive 
edges from 
 $u$ to $v$, if $A(X)_{uv}\le 0$, there are $|A(X)_{uv}|$ negative edges from 
 $u$ to $v$. 
 $\Direct{X}=(D^+(X),D^-(X))$ is such that 
\begin{align*}
 D^+(X)&=\{(i,u,v)|\ A(X)_{uv}>0, i=1\dots,A(X)_{uv}\}, \\
D^-(X)&=\{(-i,u,v)|\ A(X)_{uv}<0, i=1\dots,|A(X)_{uv}|\}.
\end{align*}
For signed directed graphs 
 $X_i=(V(X_i), \Direct{X_i}), i=1,2$, 
we define the signed directed graphs $X_1\sumg X_2$ and  
 $X_1\ddm X_2$  
by the following:
\begin{align*}
 V(X_1+X_2)&=V(X_1-X_2)=V(X_1)\cup V(X_2), \text{ and }\\
A\left(X_1+X_2\right)&=\widetilde{A}(X_1)+\widetilde{A}(X_2)
  \text{ and } 
A\left(X_1-X_2\right)=\widetilde{A}(X_1)-\widetilde{A}(X_2),
\end{align*}  
 here $\widetilde{A}(X_i)$ is the extension of $A(X_i)$ to 
$V(X_1)\cup V(X_2)$ by:
\begin{align*}
 \widetilde{A}(X_i)_{uv}=
\begin{cases}
 A(X_i)_{uv} & \text{ if }u,v\in V(X_i),\\
0 & \text{ otherwise}.
\end{cases}
\end{align*}
\end{defn} 
\begin{defn}
 Let $X=(V(X),E(X))$ be a directed graph, 
 and $B, C\subset V(X), \vv\in V(X)$.  
 We define the  signed directed graph 
 $B(\vv)\tright C$ by: 
 $V(B(\vv)\tright C)=B\cup C$, and
\begin{align*}
 A\left(B(\vv)\tright C\right)_{uv}
=\begin{cases}
 A(X)_{\vv v} & \text{ if }u\in B\text{ and }v\in C,\\
 0 & \text{ otherwise}.
\end{cases}
\end{align*}
  Let $\pi=\{V_1,\dots,V_r\}$ be a partition of vertices $V(X)$ 
 of a directed graph $X$, and 
 fix a set of vertices $\{\overline{v}_i\in V_i, i=1,\cdots,r\}$.
 Let 
\begin{align*}
 V':=\smvi{V(X)}, \ V_i':=V_i\ssetminus \{\vv_i\}. 
\end{align*}  
 We define a signed directed graph 
$\delG{X}$ 
and call it the {\em deletion graph over the partition} $\pi$, 
as follows: 
\begin{align*}
 \delG{X}&:= X\sita{V'}\ddm 
\sum_{i,j=1}^r\left(V_i'(\vv_i)\tright V_j'\right),  
\end{align*}
 here, $X\sita{V'}$ is the restriction of $X$ to $V'$.

Let  
$\Pc\in \Mat{V(X)}{\pi}{\{0,1\}}$ be the {\em  characteristic matrix}, that is, $$
\dindextwoone{\Pc}{i}{k}{j} =\delta_{ij} \text{ for }\vertex{i}{k}\in V_i.
$$  
\if01
 fix a set of deleting vetrices  $\{\overline{v}_i\in V_i, i=1,\cdots,r\}$
 and let $V':=V(X)\ssetminus \vi.$ 
\fi
For $M\in\Mat{V(X)}{V(X)}{\CC}$,
 define 
 the {\em deletion matrix over the partition} $\pi$, 
$\delG{M}$ 
$\in\Mat{V'}{V'}{\CC}$ 
 by the following: 
\begin{align*}
 \delG{M}&:=M\sita{V'\times V'} 
- \Pc\sita{V'\times \pi}\cdot M\sita{\vi\times V'}, \\
\text{ here }
 & \left(\Pc\sita{V'\times \pi}\cdot M\sita{\vi\times V'}\right)_{\vertex{i}{k}\vertex{j}{l}}=\sum_{h=1}^rP_{\vertex{i}{k}h}M_{\overline{v}_h\vertex{j}{l}}
 =M_{\vv_iv^j_l}.
\end{align*}
\end{defn}
 The following holds. 
\begin{pro}
\label{cprime}
$$\delG{A(X)}=A(\delG{X}).$$ 
\end{pro}
\text{(Proof) }\quad  
 Let $C'(X)=P\sita{V'\times \pi}\cdot A(X)\sita{\vi\times V'}.$
 Then for $\vectortwo{i}{k}\in V_i':=V_i\ssetminus \{\vv_i\}$ and 
$\vectortwo{j}{l}\in V_j'$,  we have 
 $ C'(X)_{\vectortwo{i}{k}\vectortwo{j}{l}}
= A(X)_{\overline{v}_{i}\vectortwo{j}{l}}, 
$ and 
\begin{align*}
A\left(\sum_{i,j=1}^r
\left(V_i'(\vv_i)\tright V_j'\right)\right)_{v^i_kv^j_l}
=A\left(V_i'(\vv_i)\tright V_j'\right)_{v^i_kv^j_l}= A(X)_{\vv_iv^j_l}.
\end{align*}
 Since $A\left(X\sita{V'}\right)=A(X)\sita{V'\times V'}$, 
 we have the assertion. \quad\qed
\vspace*{5mm}\\
\begin{ex}
  Consider Example \ref{join} (c) (i).
 Let $\vv_1=\vectortwo{1}{5}=$\marukk{5}{red!96!white}{.2}, 
$\vv_2=\vectortwo{2}{5}= \,$\sikakukk{5}{black}{.2}.
 Then 
\begin{align}\label{Petd}
P\sita{V'\times \pi}\cdot A(X)\sita{\vi\times V'}
\nonumber
&=
 \begin{pmatrix}
 \one_4& \zero \\
\zero &\one_4
\end{pmatrix}
\left(
\begin{array}{@{\,}cccc|cccc@{\,}}
 0&1&1& 0& 0& 0& 0& 0 \\ \hline
0& 0& 0& 0&1& 0& 0&1
\end{array}
\right) \\
\nonumber
&=
\left(
\begin{array}{@{\,}cccc|cccc@{\,}}
 \zero&\one_4&\one_4& \zero& \zero& \zero& \zero& \zero \\ \hline
\zero& \zero& \zero& \zero&\one_4& \zero& \zero&\one_4
\end{array}
\right)\\
&=
 \begin{pmatrix}
 J_{23}& \zero \\
\zero &J_{14}
\end{pmatrix},
\end{align}
 here $J_{k_1\cdots k_s}$ is the matrix such that 
 $k_1,\cdots,k_s$-column vectors are $\one$, and others are $\zero.$
 On the other hand, 
\def\oo{0.17}
\begin{align*}
 A\left(\sum_{i, j=1}^2\left(V_i'(\vv_i)\tright V_j'\right)\right)
&=
\begin{pmatrix}
 A\left(\{\marukk{\scriptsize\emph{k}}{red!96!white}{\oo}\}_{k=1}^4(\marukk{\scriptsize 5}{red!96!white}{\oo})\tright \{\marukk{\scriptsize\emph{l}}{red!96!white}{\oo}\}_{\emph{l}=1}^4\right)
&  A\left(\{\marukk{\scriptsize\emph{k}}{red!96!white}{\oo}\}_{\emph{k}=1}^4(\marukk{\scriptsize 5}{red!96!white}{\oo})
\tright \{\sikakukk{\scriptsize\emph{l}}{}{\oo}\}_{\emph{l}=1}^4\right)\\
 A\left(\{\sikakukk{\scriptsize\emph{k}}{}{\oo}\}_{\emph{k}=1}^4(\sikakukk{\scriptsize 5}{}{\oo})
\tright \{\marukk{\scriptsize\emph{l}}{red!96!white}{\oo}\}_{\emph{l}=1}^4\right)
&   A\left(\{\sikakukk{\scriptsize\emph{k}}{}{\oo}\}_{\emph{k}=1}^4(\sikakukk{\scriptsize 5}{}{\oo})
\tright \{\sikakukk{\scriptsize\emph{l}}{}{\oo}\}_{\emph{l}=1}^4\right)
\end{pmatrix}\\
&=
\begin{pmatrix}
 J_{23}& \zero \\
\zero &J_{14}
\end{pmatrix}.
\end{align*}
 Next, consider Example \ref{join} (c) (ii).
 Let $\vv_1=\vertex{1}{1}=\marukkrotate{1}{}{.2}\ ,\ 
\vv_2=\vertex{2}{3}=\maruk{3}{red!96!white},
\vv_3=\vertex{3}{6}=\sikakuk{6}{black}$. 

 Then 
\begin{align*}
 P\sita{V'\times \pi}\cdot A(X)\sita{\vi\times V'}
&=\left(
\begin{array}{@{\,}c|c|c@{\,}}
 \zero&\one_2& \zero \\ \hline
\zero& \zero&\one_5
\end{array}
\right)
\left(
\begin{array}{@{\,}cc|cccccc@{\,}}
 1&1& 0& 0& 0& 0& 0 \\ 
 0&0& 0& 0& 0& 0& 1 \\ 
 0&0& 0& 1& 0& 1& 0 
\end{array}
\right)\\
&=\left(
\begin{array}{@{\,}cc|ccccc@{\,}}
 \zero& \zero& \zero& \zero& \zero& \zero&\one_2 \\ 
 \zero& \zero& \zero& \one_5& \zero& \one_5&\zero 
\end{array}
\right)=
\begin{pmatrix}
 \zero&J_{5} \\
\zero &J_{24}
\end{pmatrix}.
\end{align*}
 On the other hand, 
\begin{align*}
& A\left(\sum_{i, j=1}^3\left(V_i'(\vv_i)\tright V_j'\right)\right)\\
&=
\begin{pmatrix}%
 A\left(\{\marukk{\scriptsize 1}{red!96!white}{\oo},\marukk{\scriptsize 2}{red!96!white}{\oo}\}(\marukk{\scriptsize 3}{red!96!white}{\oo})
\tright \{\marukk{\scriptsize 1}{red!96!white}{\oo},\marukk{\scriptsize 2}{red!96!white}{\oo}\}\right)
&  A\left(\{\marukk{\scriptsize 1}{red!96!white}{\oo},\marukk{\scriptsize 2}{red!96!white}{\oo}\}(\marukk{\scriptsize 3}{red!96!white}{\oo})
\tright \{\sikakukk{\scriptsize 1}{}{\oo},\dots,\sikakukk{\scriptsize 5}{}{\oo}\}\right)\\
A\left(\{\sikakukk{\scriptsize 1}{}{\oo},\dots,\sikakukk{\scriptsize 5}{}{\oo}\}(\sikakukk{\scriptsize 6}{}{\oo})
\tright \{\marukk{\scriptsize 1}{red!96!white}{\oo},\marukk{\scriptsize 2}{red!96!white}{\oo}\}\right)
&A\left(\{\sikakukk{\scriptsize 1}{}{\oo},\dots,\sikakukk{\scriptsize 5}{}{\oo}\}(\sikakukk{\scriptsize 6}{}{\oo})
\tright \{\sikakukk{\scriptsize 1}{}{\oo},\dots,\sikakukk{\scriptsize 5}{}{\oo}\}\right)
\end{pmatrix}\\
&=
\begin{pmatrix}
\zero& J_{5} \\
\zero &J_{24}
\end{pmatrix}. 
\end{align*} 
 Hence, $\delG{A(X)}=A(\delG{X}).$  \qed
\end{ex}
\fi
\section{The Similarity Transformation}
\begin{lem}\label{MP}
 Let $(M,\pi)$ be an equitable pair, here $M\in \Mat{V(X)}{V(X)}{\mathbb{C}}$ 
and  $\pi=\{V_1,\dots,V_r\}$ 
 is  a partition of vertex set $V(X)$ of a graph $X$.  
Then 
$$M\Pc=\Pc\cdot M\quotient{\pi}.$$ 
\end{lem}
%
\begin{align*}
\mathrm{(Proof)} \qquad (MP)_{\vt{i}{k}j}
&=\sum_{\vt{h}{l}\in V}
\dindex{M}{i}{k}{h}{l}
\dindextwoone{P}{h}{l}{j}
=\sum_{l}\dindex{M}{i}{k}{j}{l}=(M\quotient{\pi})_{ij},
\text{ and }\\ 
\dindextwoone{(P\cdot M\quotient{\pi})}{i}{k}{j}  
&=\sum_{h=1}^r\dindextwoone{P}{i}{k}{h}(M\quotient{\pi})_{hj}
=(M\quotient{\pi})_{ij}.\quad \qed 
\end{align*}
\begin{defn}
 Let $\pi=\{V_1,\dots,V_r\}$ be a partition of vertices $V(X)$ 
 of a graph $X$.
 Fix $\{\overline{v}_i\in V_i, i=1,\cdots,r\}$.
 Define 
$$Q=Q(\vi)\in\Mat{V(X)}{V'}{\{0,1\}},$$ 
 here  $V'=\smvi{V(X)}$, by 
 the following:
$$  
\dindex{Q}{i}{k}{j}{l}
:=\dindex{\delta}{i}{k}{j}{l}.
$$
 For instance, in Example \ref{join} (a), letting $\vv_i=v^i_2, i=1,2,$ 
\begin{align*}
 P=
\begin{pmatrix}
 1& 0\\1& 0\\0& 1\\0& 1
\end{pmatrix}, \text{ and }
 Q=
\begin{pmatrix}
 1& 0\\0& 0\\0& 1\\0& 0
\end{pmatrix}.
\end{align*}
\end{defn}
\begin{lem}
 The   column vectors of $(\Pc,Q)\in \Mat{V(X)}{(\pi\times V')}{\{0,1\}}$ 
 form a basis of $\CC^{V(X)}$, that is, $(\Pc,Q)$ is invertible. 
\end{lem}
(Proof)\quad Assume
\begin{align*}
 \sum_{i=1}^r\lambda_iP_{\cdot i}+  \sum_{\vectortwo{j}{l}\in V_j',j=1,\dots,r}
\mu_{\vectortwo{j}{l}}Q_{\cdot\vectortwo{j}{l}}=0,
\end{align*}
 here $P_{\cdot i}$ (resp. $Q_{\cdot\vectortwo{j}{l}}$) 
 is the $i$-th (resp. $\vectortwo{j}{l}$-th) column vector of $P$ 
 (resp. $Q$). 
 For $s\in \{1,\dots,r\}$, 
 $P_{\vv_si}=\delta_{si}$, and since 
$v^j_l\in V_j'\subset V'=V(X)\ssetminus \{\vv_i\}_{i=1}^r$,
 we have $v^j_l\neq \vv_s$, implying $Q_{\vv_sv^j_l}=0.$ 
 So, considering the entry corresponding to $\vv_s$ we have $\lambda_s=0, s=1,\dots,r.$ 
 And we have 
 $\mu_{\vectortwo{j}{l}}=0, \vectortwo{j}{l}\in V', j=1,\dots, r.$ \medskip\qed 
\newcommand{\PP}{\overline{P}}

 The following is our main theorem. 
\begin{thm}\label{zeta}
 Let $(M,\pi)$ be an equitable pair, here $M\in \Mat{V(X)}{V(X)}{\mathbb{C}}$ 
and  $\pi=\{V_1,\dots,V_r\}$ 
 is  a partition of vertex set $V(X)$ of a graph $X$.  
 Then letting $\PP=(P,Q)$ we have 
\begin{align*}
 \PP^{-1}M \PP= \begin{pmatrix}
 M\quotient{\pi} & C\\
  \zero & \delG{M}
\end{pmatrix}, C=M\sita{\vi\times V'}\,,
\end{align*}
 hence, 
\begin{align*}
 \phi(M,x)&=\phi(M\quotient{\pi},x)\cdot \phi(\delG{M},x), 
\end{align*}%
 here $\delG{M}$ is defined for a fixed 
$\vi, \vv_i\in V_i$,  
 and $\phi(\delG{M},x)$ does not depend on 
$\vi$.
In particular, if $\pi$ is equitable, then 
 \begin{align*}
\charp{A(X)}= \charp{A(X\quotient{\pi})}\cdot
 \charp{A\left(\delG{X}\right)}.
 \end{align*}
\end{thm}
(Proof)\quad  
 We show that 
\begin{align*}
 MQ&= P C+Q\cdot \delG{M}.
\end{align*}
 For 
$v^i_k\in V_i, 
\mbox{\footnotesize $\vectortwo{j}{l}$}\in V_j'=V_j\ssetminus \{\vv_j\},$ 
\begin{align}\label{QQ}
 \dindex{(MQ)}{i}{k}{j}{l}
&=\sum_{w\in V}\dindextwoone{M}{i}{k}{w}
\dindexonetwo{Q}{w}{j}{l}= \dindex{M}{i}{k}{j}{l}, \text{ and }
\end{align}
\begin{align*}
 \dindex{(P C+Q\cdot \delG{M})}{i}{k}{j}{l}
&=\sum_{h=1}^{r}\dindextwoone{\Pc}{i}{k}{h}
M_{\overline{v}_hv^j_l} 
+ \sum_{w\in V'}Q_{v^i_kw}(\delG{M})_{wv^j_l}\\
&=
\begin{cases}
M_{\overline{v}_iv^j_l} 
+ (\delG{M})_{v^i_kv^j_l}  & \text{ if } v^i_k\neq \overline{v}_i\\
M_{\overline{v}_iv^j_l}  & \text{ if } v^i_k = \overline{v}_i.
\end{cases}
\\
 \text{ If } v^i_k\neq \overline{v}_i,\qquad
 (\delG{M})_{v^i_kv^j_l}
&=\dindex{\left(M\sita{V'\times V'}\right)}{i}{k}{j}{l}
-\sum_{h=1}^{r}\dindextwoone{\left(\Pc\sita{V'\times \pi}\right)}{i}{k}{h}
\dindexonetwo{M}{\overline{v}_h}{j}{l}\\
&=M_{v^i_kv^j_l}-M_{\overline{v}_iv^j_l}.
\end{align*}
\vspace{-2mm}
So that we have 
$\dindex{(P C+Q\cdot \delG{M})}{i}{k}{j}{l}=\dindex{M}{i}{k}{j}{l}
=\dindex{(MQ)}{i}{k}{j}{l}$ (by (\ref{QQ})). \quad\qed
\begin{ex}
Let $X, \pi$ be as in Example \ref{join}. 
\\
(a) Delete $\{\vv_1,\vv_2\}=\{\marukk{2}{}{.18},\sikakukk{2}{}{.18}\}$. 
\def\oo{0.15}

$X\ssetminus \pi=X\sita{\,\marukk{1}{}{\oo},\sikakukk{1}{}{\oo}}
\ 
-
\Bigl(
\left(\{\marukk{1}{}{\oo}\}\left(\marukk{2}{}{\oo}\right)\tright
\marukk{1}{}{\oo}\right)
+\left(\{\marukk{1}{}{\oo}\}\left(\marukk{2}{}{\oo}\right)\tright\sikakukk{1}{}{\oo}\right)$

\hspace{55mm}$+ \left(\{\sikakukk{1}{}{\oo}\}\left(\sikakukk{2}{}{\oo}\right)\tright\marukk{1}{}{\oo}\right)
+\left(\{\sikakukk{1}{}{\oo}\}\left(\sikakukk{2}{}{\oo}\right)\tright\sikakukk{1}{}{\oo}\right)
\Bigr)
$ 
\vspace*{-8mm}

\hspace*{11mm}
\def\ooo{.85}
$=$ 
\CtwoK{1}{1}{0.2}{}{-1.6,0}{1,1}{0.6}{\ooo} 
\hspace*{2mm}
$-$ \hspace*{-3mm}
\KakkoLL{0}{0.8}{-.3,-0.1}{.5,1.8}{\ooo}
\hspace*{-3mm}
\ConeMs{1}{1}{0.2}{}{-1,-0.1}{1.8,1.8}{0.5}{\ooo}
\hspace*{-10mm}
 $+$ \hspace*{-6mm}
\CtwoNKs{1}{1}{0.2}{}{-2,-0.1}{1.8,1.8}{0.6}{\ooo}
\hspace*{-5mm} 
$+$ 
\hspace*{-6mm} 
\CtwoNKs{1}{1}{0.2}{}{-2,-0.1}{1.8,1.8}{0.5}{\ooo}
\hspace*{-5mm} 
 $+$ 
\hspace*{-6mm} 
\ConeKs{1}{1}{0.2}{}{-1,-0.1}{1.8,1.8}{0.4}{\ooo}
\hspace*{-17mm}
\KakkoRs{0}{0.8}{-.3,-0.1}{.5,1.8}{\ooo}
\vspace*{.1mm}

\hspace*{11mm}
$=$ 
\hspace*{-4mm} 
\CtwoLoopKs{1}{1}{0.2}{}{-2,-0.1}{1.8,1.8}{0.6}{\ooo}\hspace*{-1.5mm}, 
\hspace*{-1.1mm}
 here we denoted the negative edeges by the red lines.
\vspace*{3mm}

%
 Hence, $A(X\ssetminus\pi)=
\begin{pmatrix}
 -1& 1\\
 1 & -1
\end{pmatrix}.$
 By our Theorem, we have 
 $\charp{A(X)}\\ =
\begin{vmatrix}
 x-1& -1\\
 -1 & x-1
\end{vmatrix}\cdot
\begin{vmatrix}
 x+1& -1\\
 -1 & x+1
\end{vmatrix}=((x-1)^2-1)\cdot((x+1)^2-1)=x^2(x+2)(x-2).
$\\
\if01
 By the direct computation, we have 
$  \charp{A(X)}=
\begin{vmatrix}
 xI_2-A_1& -I_2& \\
-I_2&xI_2-A_1
\end{vmatrix}\\
=|(xI_2-A_1)^2-I_2|
=|(x+1)I_2-A_1|\cdot|(x-1)I_2-A_1|
=((x+1)^2-1)\cdot((x-1)^2-1).
$
\fi
\vspace{0mm}

(b) Delete $\{\vv_1,\vv_2\}=\{\marukk{3}{}{.18},\sikakukk{2}{}{.18}\}$. 
\def\oo{0.15}
\vspace*{2mm}
 
$X\ssetminus \pi$
\hspace*{0mm}
$=X\sita{\,\marukk{1}{}{\oo},\marukk{2}{}{\oo},\sikakukk{1}{}{\oo}}
\ -\Bigl(
\left(\{\marukk{1}{}{\oo},\marukk{2}{}{\oo}\}\left(\marukk{3}{}{\oo}\right)\tright
\{\marukk{1}{}{\oo},\marukk{2}{}{\oo}\}\right)
+\left(\{\marukk{1}{}{\oo},\marukk{2}{}{\oo}\}\left(\marukk{3}{}{\oo}\right)
\tright\sikakukk{1}{}{\oo}\right)$
\vspace*{2mm}\\
\hspace{5mm}
$+\left(\{\sikakukk{1}{}{\oo}\}\left(\sikakukk{2}{}{\oo}\right)
\tright\{\marukk{1}{}{\oo},\marukk{2}{}{\oo}\}\right)
+\left(\{\sikakukk{1}{}{\oo}\}\left(\sikakukk{2}{}{\oo}\right)
\tright\{\sikakukk{1}{}{\oo}\}\right)
\Bigr)$

\vspace*{12mm}

\def\rai{15mm}
\hspace*{-4mm}\raisebox{\rai}{$=$} 
%
\begin{tikzpicture}[scale=0.7,>=stealth]
\useasboundingbox (-1,-1) rectangle (1,1);
\LKakko{0.8}{0.5}
\Dfivethree
\RKakko{0.8}{0.5}{12}
\end{tikzpicture}
 \hspace*{11mm} \raisebox{\rai}{$-$} 
\def\r{0.25}\Mul\r{0.65}\acomo\def\c{0.9}\Div\r{1.2}\rr
\def\hh{0.8}\def\ll{1.5}\def\he{0.8}\Add\ll{0.2}\lll
\Div\he{3}\hehe
\begin{tikzpicture}[scale=0.7,>=stealth]
\useasboundingbox (-1,-1) rectangle (.5,.5);
\LKakko{0.8}{0.5}
\hspace*{3mm}
\LKakko{0.6}{0.5}
\def\R{0.2}
\Node{0}{2}{1}{\R}
\Node{0}{0}{2}{\R}
\EdgeOneUU
\LoopUDtwo{0}{2}{}{}
\end{tikzpicture}
\hspace{-20mm}
\begin{tikzpicture}[scale=0.7,>=stealth]
\useasboundingbox (-1,-1) rectangle (.5,.5);
\RKakko{0.6}{0.5}{12}
\end{tikzpicture}
\hspace{15mm}\raisebox{\rai}{$+$} 
\begin{tikzpicture}[scale=0.8,>=stealth]
\useasboundingbox (-1,-1) rectangle (1,1);
\LKakko{0.6}{0.5}
\def\R{0.2}
\Node{0}{2}{1}{\R}
\Node{0}{0}{2}{\R}
\NodeRec{2}{2}{1}{\R}
\end{tikzpicture}
\hspace*{-15mm}
\begin{tikzpicture}[scale=0.8,>=stealth]
\useasboundingbox (-1,-1) rectangle (1,1);
\RKakko{0.6}{0.5}{12}
\end{tikzpicture}
\hspace{11mm}\raisebox{\rai}{$+$}  
\begin{tikzpicture}[scale=0.8,>=stealth]
\useasboundingbox (-1,-1) rectangle (1,1);
\LKakko{0.5}{0.5}
\Sub\r{0.1}\rrr\def\R{0.2}
\Node{0}{2}{1}{\R}
\Node{0}{0}{2}{\R}
\NodeRec{2}{2}{1}{\R}
\Sub{2}\rr\Reclu
\def\ruldac{.85}
\draw[very thick,->] (\Reclu,\Reclu)--(1,1);
\draw[very thick,->] (1.1,1.1)--(\ruldac,\ruldac);
\draw[very thick] (1,1)--(\rrr,\rrr);
\end{tikzpicture}
\hspace*{-15mm}
\begin{tikzpicture}[scale=0.8,>=stealth]
\useasboundingbox (-1,-1) rectangle (1,1);
\RKakko{0.5}{0.5}{12}
\end{tikzpicture}
\hspace*{11mm}\raisebox{15mm}{$+$}\hspace*{-11mm} 
\begin{tikzpicture}[scale=.8,>=stealth]
\useasboundingbox (-2,-2) rectangle (0,0);
 \ConeKK{1}{2}{0.2}{}{0.4} 
\end{tikzpicture}
\hspace*{-7mm}
\begin{tikzpicture}[scale=.8,>=stealth]
\useasboundingbox (-1,-1) rectangle (1,1);
\RKakko{0.6}{0.5}{1}
\end{tikzpicture}
\vspace*{8mm}

\hspace*{-4mm}\raisebox{18mm}{$=$} 
\begin{tikzpicture}[scale=0.8,>=stealth]
\useasboundingbox (-1,-1) rectangle (1,1);
\LKakko{0.8}{0.5}
\def\R{0.2}
\Node{0}{2}{1}{\R}
\Node{0}{0}{2}{\R}
\NodeRec{2}{2}{1}{\R}
\Edge{0}{2}{0}{0}{\R}{red}
\EdgeCRecrll{0}{2}{2}{2}{\R}
\Sub{2}\rr\Reclu
\Sub\r{0.1}\rrr\def\ruldac{.85}
\draw[very thick,->,red] (\Reclu,\Reclu)--(1,1);
\draw[very thick,->,red] (1.1,1.1)--(\ruldac,\ruldac);
\draw[very thick,red] (1,1)--(\rrr,\rrr);
\LoopU{0}{2}{1.2}{red}
\LoopU{2}{2}{1.2}{red}
\LoopD{0}{0}{1.2}{}
\end{tikzpicture}
\hspace*{-4mm}
\begin{tikzpicture}[scale=.8,>=stealth]
\useasboundingbox (-1,-1) rectangle (1,1);
\RKakko{0.8}{0.5}{1}
\end{tikzpicture}
\raisebox{13mm}{$.$ }
\vspace*{-30mm}

\begin{align*}
\hspace{21mm} \text{Hence, } A(X\ssetminus\pi)=
\begin{pmatrix}
 -1& 1& 1 \\
 -1 &1& 0\\
2&-2&-1
\end{pmatrix},
\text{ implying } 
\end{align*}
\begin{align*}
 \charp{A(X)}&=
\begin{vmatrix}
 x-3& -1\\
 -2 & x-1
\end{vmatrix}\cdot
\begin{vmatrix}
 x+1& -1& -1\\
 1 & x-1&0\\
-2&2&x+1
\end{vmatrix}=(x^2-4x+1)\cdot x(x-1)(x+2).
\end{align*}
(c) (i) 
 By (\ref{Petd}) we have
\begin{align*}
 A(X\ssetminus \pi_1)&=
A(X\sita{V'})- P\sita{V'\times \pi_1}\cdot A(X)\sita{\{\vv_i\}_{i=1}^2\times V'}\\
&=
 \begin{pmatrix}
 A_1'& I_4& \\
I_4&A_2'
\end{pmatrix}
-\begin{pmatrix} J_{23}& \zero \\
\zero &J_{14}\end{pmatrix}
=\begin{pmatrix}
 A_1& I_4& \\
I_4&A_2
\end{pmatrix},
\end{align*} 
\begin{align*}
\text{here } A_1&=A_1'-J_{23}= 
\begin{pmatrix}
 0&0&1&1 \\
0&0&0&1\\
1&0&0&0\\
1&1&0&0
\end{pmatrix}-J_{23}=
\begin{pmatrix}
 0&-1&0&1 \\
0&-1&-1&1\\
1&-1&-1&0\\
1&0&-1&0\\
\end{pmatrix},
\end{align*}
\begin{align*}
A_2&=A_2'-J_{14}= 
\begin{pmatrix}
 0&1&0&0 \\
1&0&1&0\\
0&1&0&1\\
0&0&1&0
\end{pmatrix}-J_{14}=
\begin{pmatrix}
 -1&1&0&-1 \\
0&0&1&-1\\
-1&1&0&0\\
-1&0&1&-1
\end{pmatrix}.
\end{align*}
\begin{align*}
\text{ So, } \charp{A(X\ssetminus \pi_1)}
&=\begin{vmatrix}
 xI_4-A_1&-I_4 \\
-I_4&xI_4-A_2
\end{vmatrix}
=|(xI_4-A_1)(xI_4-A_2)-I_4|\\
&=|x^2I_4-x(A_1+A_2)+A_1A_2-I_4|.
\end{align*}
 Since $A_1+A_2=A_1A_2=-I_4$, we have 
$\charp{A(\delG{X})}=|(x^2+x-2)I_4|=(x^2+x-2)^4.$ 
 By $\charp{A(X\quotient{\pi_1})}=(x-3)(x-1)$ we have  
\begin{align*}
 \charp{A(X)}=(x-3)(x-1)(x^2+x-2)^4=(x-1)^5(x+2)^4(x-3).
\end{align*} 
 To get the characteristic polynomial of the adjacency matrix 
 of the Petersen graph by the another method, see  [GR01, section 9.1 and exercise 8.9].
\end{ex}
\if00
\section{Applications to the Ihara-Bartholdi zeta functions}
 In this section, we assume $X=(V(X), E(X))$ is an undirected graph,   
 that is, $E(X)$ is a set of symmetric directed edges. 
 As corollaries of our Theorem, we have the decomposition 
 of the characteristic polynomial of Laplacian matrix, 
 and  the reciprocal of the Ihara-Bartholdi zeta functions 
 of equitably partitioned graphs.

 Zeta functions of a graph are defined as follows.
 For an edge $e\in E(X)$, 
 we denote by $o(e)$  (resp. $t(e)$) the origin  (resp. terminus) of $e$.
A closed path in $X$ is an sequence of edges $C=(e_1,\dots,e_k)$ 
 with $t(e_i)=o(e_{i+1})$ for $i\in \ZZ/k\ZZ$. 
 We denote by $k=|C|$, the length of $C$, 
and by $cbc(C)=\#\{i\in \{1,\dots,k\}|\ e_{i+1}=\overline{e_i}\}$, 
 the cyclic bump count of $C$. 
 A cycle $[C]$ is the equivalence class of a closed path $C$ under cyclic 
 permutation of its edges 
  (that is, $(e_1,\dots,e_k)\sim (e_2,\dots,e_k, e_{1})$). 
 A cycle  is {\em prime} if none of its representatives can be written 
 as  $C^k$ for some $k\ge 2$.
 We denote by $\mathcal{C}$ the set of prime cycles. 
 Bartholdi zeta function is defined by 
$$Z_X(u,t)=\prod _{[C]\in \mathcal{C}}\frac{1}{1-u^{cbc(C)}t^{|C|}},$$

 $Z_X(0,t)=Z_X(t)$ is the Ihara zeta function defined 
by Ihara [Iha66] in which 
 he considered  a zeta function of a regular graph and 
 gave its reciprocal as a polynomial. 
 It was generalized to 
 general graphs by [Bas92] (see also [Ser80], [Has89], [Has90], 
 [Sun96a], [Sun96b],  [ST96], [FZ99], [KS00]). 
 Bartholdi generalized Bass's Theorem as the following. 
\begin{thm}\label{Bar}[Bar99]
 Let $X$ be a connected graph with $n$ vertices and 
 $m$ (non-oriented) edges. Then the reciprocal of the  Bartholdi-Ihara 
 zeta function of $X$ is given by  
\begin{align*}
 & Z_X(u,t)^{-1}\\
&=\left(1-(1-u)^2t^2\right)^{m-n}
 \det \left(I_n-tA(X)+(1-u)(D(X)-(1-u)I_n)t^2\right),
\end{align*}
 where $A(X)$ is the adjacency matrix of $X$, and $D(X)=(d_i\delta_{ij})$ 
 is the degree matrix which is diagonal with $d_i=\mathrm{deg}v_i, 
 V(X)=\{v_1,\dots,v_n\}$. 
\end{thm}
For $u, t\in\CC,$ let 
\begin{align}\label{zetap}
 \Dz&=s_1I_n+s_2 D(X)\in\Mat{V(X)}{V(X)}{\mathbb{C}} 
\end{align}
 be the diagonal matrix, here $s_1 =1-(1-u)^2t^2, \ s_2=(1-u)t^2.$
Then we have $$Z_X(u,t)^{-1}=s_1^{m-n}\det (-t A(X)+\Dz).$$
 As a Corollary of  Theorem  \ref{zeta}, 
 we have the following. 
\begin{cor}\label{Let}
    Let $\pi=\{V_1,\dots,V_r\}$ be an equitable partition of the 
 vertex set $V(X)$ 
 of a graph $X$ and fix a set of deleting vertices 
 $\{\overline{v}_i\in V_i, i=1,\cdots,r\}$. 
 Let $D$ be a diagonal matrix such that 
  $D\sita{V_i\times V_i}=d_iI_i$ 
 (for instance, $D={\large O}, D(X)$ or $D^Z(X)$).   
 Then we have 
\begin{align*}
 \phi\left(\alpha A(X)+D,x\right)
&= \phi\left(\alpha A(X\quotient{\pi})+(\delta_{ij}d_i),x\right)
\cdot\phi\left(\alpha A(\delG{X})+D\sita{V'\times V'},x\right).
\end{align*}
In particular, 
 \begin{align*}
 (a)\quad 
& \charp{L(X)}=\charp{-A(X\quotient{\pi})+D(X\quotient{\pi})}\cdot
\charp{-\CX+\rVd{D(X)}},\\
& \charp{Q(X)}=
\charp{A(X\quotient{\pi})+D(X\quotient{\pi})}\cdot
\charp{A(\delG{X})+\rVd{D(X)}}.\\
(b)\quad  & Z_X(u,t)^{-1}\\
&=s_1^{m-n}
\det \left(-tA(X\quotient{\pi})+D^Z(X\quotient{\pi})\right)\cdot
\det\left(-tA(\delG{X})+\rVd{D^Z(X)}\right),
\end{align*}  
 $\text{here }
  s_1 =1-(1-u)^2t^2.$ 
\end{cor}
(\text{Proof}) \quad  \\ 
 Let $M=\alpha A(X)+D$. 
Then $M\quotient{\pi}=\alpha A(X\quotient{\pi})+(\delta_{ij}d_i)$.
For $\vt{j}{l}\in V_j':=V_j\ssetminus \{\vv_j\},$   
$\dindexonetwo{\left(D\sita{\vi\times V'}\right)}{\vv_i}{j}{l}
=\dindexonetwo{D}{\vv_i}{j}{l}=0.   
$ So, we have 
$
M\sita{\vi\times V'}=\alpha A(X)\sita{\vi\times V'}.
$
 Hence, 
\begin{align*}
 \delG{M}
&= \alpha A(X\sita{V'\times V'})+D\sita{V'\times V'}
  -P\sita{V'\times\pi}\cdot \alpha A(X)\sita{\vi\times V'}\\
&= \alpha (\delG{A(X)})+ D\sita{V'\times V'}
= \alpha A(\delG{X})+ D\sita{V'\times V'} 
\ (\text{by Proposition \ref{cprime}}).
\end{align*}
 So, we have the first assertion.
By letting 
(a)  $D=D(X)$, (i) $\alpha=-1,$ 
(ii) $\alpha=1,$ (b) $D=D^Z(X), \alpha=-t$,   
 we have the other assertions. \quad  \qed
\vspace{4mm}\\
\if01
\textbf{(Proof of Corollary \ref{Let})} \  
 Let (a) $\alpha=1, D=0$, 
(b) $D=D(X)$, (i) $\alpha=-1$, (ii) $\alpha=1$, 
 and (c) $\alpha=-t, D=D^Z(X)$, we have the assertions. \qed
\vspace{4mm}
\\ 
\fi
\begin{rem}
 Note that $\rVd{D(X)}\neq D(X\sita{V'})$ 
 in general. Let $X$ be a cycle $C_4$ in Example \ref{join} (a). 
Then  $V'=\{\marukk{1}{}{.18},\sikakukk{1}{}{.18}\}, D(X)=2I_4$, and 
$D(X)\sita{V'\times V'}
=2I_2$, but $D(X\sita{V'})=D(C_2)=I_2.$ 
\if01
$V(X)=\{0,1,2,3\}, i\sim j$ iff $i-j\equiv \pm 1$ mod  $4.$
 Then $D(X)=2I_4.$ Let $\pi=\{\{0,1\},\{2,3\}\}$, 
and $v_1=1,v_2=2$. 
\fi
 \end{rem}
\section{Generalized join graphs} 
 In this section, graphs considered are assumed to be simple. 
 We consider the case when  $X$ is the generalized join (composition)
  of $X_1,\dots,X_r$ (see Example \ref{join} (d)). 
 The following Corollary of Theorem \ref{zeta} includes  results of [Sch74, Theorem 7] and 
 [CC17, Theorem 3.1].  
\begin{cor}\label{jjoin}
\quad  Let $X=H[X_1,\dots,X_r]$ be a generalized join (composition) 
with each $X_i$ 
being $k_i$-regular. 
 Let $D\in\Mat{V(X)}{V(X)}{\mathbb{C}}$ be a diagonal matrix such that 
$D\sita{V(X_i)\times V(X_i)}=d_iI_i.$ 
 Then, letting $\pi=\{V(X_1),\dots,V(X_r)\}$, for any $\alpha\in\mathbb{C},$
\begin{align}\label{Xover} 
\charp{\alpha A(X)+D}%
&=\charp{\alpha A(X\quotient{\pi})+(\delta_{ij}d_i)}
\cdot\prod_{i=1}^r\frac{\charp{\alpha A(X_i)+d_iI_i}}{x-\alpha k_i-d_i}\\
\label{H}&=\det\left(-\alpha A(H)+\DD(x)\right)
\cdot\prod_{i=1}^r\frac{n_i\charp{\alpha A(X_i)+d_iI_i}}
{x- \alpha k_i-d_i},\\
\nonumber
\hspace{-5mm}\text{ here }\DD(x)_{ij}&=\delta_{ij}\frac{1}{n_i}(x- \alpha  k_i-d_i),  n_i=\#V(X_i).    
\end{align} 
In particular,
\begin{align*}
\mathrm{(a)}\quad  \charp{A(X)} 
&=\charp{A(X\quotient{\pi})} 
\cdot\prod_{i=1}^r\frac{\charp{A(X_i)}}{x-k_i}\hspace{40mm}\\
&= \det \left(-A(H)+\DD(x)\right)
\cdot\prod_{i=1}^r\frac{n_i\charp{A(X_i)}}{x-k_i},\\
\text{here }\DD(x)_{ij}&=\delta_{ij}\frac{1}{n_i}(x-k_i), \text{  and }
\end{align*}
\begin{align}\label{gjoin}
\nonumber \mathrm{(b)}\quad   Z(u,t)^{-1}&\\
=s_1^{m-n}&
\det\left(-tA(X/\pi)+ (\delta_{ij}d_i^Z(X))\right)\cdot 
\prod_{i=1}^r\frac{\det (-tA(X_i)+d_i^Z(X)I_i)}{-tk_i+\dzi}\\
\nonumber =s_1^{m-n}&
\det(-tA(H)+\DD)
\prod_{i=1}^r\frac{n_i\det (-tA(X_i)+d_i^Z(X)I_i)}{-tk_i+\dzi}, 
\end{align}
 here 
$\DD\in\Mat{r}{r}{\mathbb{C}}$ is a diagonal matrix 
such that 
 $\DD_{ii}= \frac{1}{n_i}(-tk_i+\ddz_i(X))\in\mathbb{C}, 
\ddz_i(X)
=s_1+s_2d_i(X)=s_1+s_2(k_i+\sum_{j\neq i}A(H)_{ij}n_j), s_1={1-(1-u)^2t^2,}$ 
$s_2=(1-u)t^2, d_i(X)=\deg v, v\in V(X_i). $
\end{cor}
\text{(Proof)}\quad 
Let $M=\alpha A(X)+D$. 

(i)  First we prove that 
 \begin{align}\label{CC}
 \charp{\delG{M}} 
=\prod_{i=1}^r\frac{\charp{\alpha A(X_i)+d_iI_i}}{x-\alpha k_i-d_i}, 
\end{align}
 which implies (\ref{Xover}) by Theorem \ref{zeta} and (\ref{equitable}). 
 By Example \ref{join} (d) and Remark \ref{equi}, 
$(M,\pi)$ is an equitable pair. 
 Let $V_i':=V(X_i)\ssetminus \{\vv_i\}, i=i,\dots, r.$ 
 As  the notation in Theorem \ref{zeta}, 
 for $i\neq j$, 
$\vt{i}{k}\in V_i',\vt{j}{l}\in V_j',$
 since 
$ \left(P\sita{V'\times\pi}\cdot A(X)\sita{\vi\times V'}\right)_{\vertex{i}{k}\vertex{j}{l}}
=A(X)_{\overline{v}_{i}\vertex{j}{l}}=A(H)_{ij},
$ we have
$A(\delG{X})_{\vertex{i}{k}\vertex{j}{l}}
=A(X)_{\vertex{i}{k}\vertex{j}{l}}
 - A(H)_{ij}=0,$ which implies 
$(\delG{M})_{\vertex{i}{k}\vertex{j}{l}}
=\alpha A(\delG{X})_{\vertex{i}{k}\vertex{j}{l}}
+D_{\vertex{i}{k}\vertex{j}{l}}=0.$

 Hence, 
$\delG{M}= \begin{pmatrix}
 C_{1}  && $\raisebox{-2ex}{\hspace{-8mm} \huge O}$ \\
&\ddots&\\
 $\raisebox{2ex}{\hspace{2mm}\huge O}$ && C_r
\end{pmatrix}, 
$ 
here $C_i=\left(\delG{M}\right)\sita{V_i'\times V_i'}$.
So, 
\begin{align}\label{C}
 \charp{\delG{M}}=\prod_{i=1}^r\charp{C_i}.
\end{align}
We prove that 
\begin{align}\label{one}
 \charp{C_i} =\frac{\charp{M_i}}{x-\alpha k_i-d_i}, 
 \text{ here } M_i=M\sita{V_i'\times V_i'}
=\alpha A(X_i)+d_iI_i, 
\end{align}
 which implies (\ref{CC}).
For a subgraph $X_i$, let $\overline{\pi}_i=\{V(X_i)\}$ be 
 the trivial partition. 
 Since 
 $M_i\one = (\alpha k_i+d_i)\one$, 
 we have 
$M_i\quotient{\overline{\pi}_i}=\alpha k_i+d_i\in\mathbb{C}$.  
 So, by Theorem \ref{zeta} we have 
\begin{align*}
 \charp{M_i}  &= (x-\alpha k_i-d_i)\charp{M_i\backslash\overline{\pi}_i}. 
\end{align*}
 We show $M_i\backslash\overline{\pi}_i=(\delG{M})\sita{V_i'\times V_i'}=C_i$, which implies (\ref{one}).
Since 
\begin{align*}
 \left((\delG{M})\sita{V_i'\times V_i'}\right)_{\vertex{i}{k}\vertex{i}{l}}
&= \left(\delG{M}\right)_{\vertex{i}{k}\vertex{i}{l}}
= \left(M\sita{V'\times V'}-P\sita{V'\times \pi}\cdot M\sita{\vi\times V'}\right)_{\vertex{i}{k}\vertex{i}{l}}\\
&= M_{\vertex{i}{k}\vertex{i}{l}}-M_{\overline{v}_i\vertex{i}{l}},
\end{align*}
 and 
\begin{align*}
\left(M_i\backslash\overline{\pi}_i\right)_{\vertex{i}{k}\vertex{i}{l}}
&=\left(M_i\sita{V'_i\times V'_i}-P\sita{V'_i\times \overline{\pi}_i}\cdot 
M_i\sita{\{\overline{v}_i\}\times V'_i}\right)_{\vertex{i}{k}\vertex{i}{l}}
= M_{\vertex{i}{k}\vertex{i}{l}}-M_{\overline{v}_i\vertex{i}{l}}.
\end{align*}
 So, we have the assertion. 

(ii) Next we prove that 
\begin{align}\label{M}
 \charp{M\quotient{\pi}} =n_1\cdots n_r\det\left(-\alpha A(H)+\DD(x)\right),
\end{align}
 which implies (\ref{H}). 
 By (\ref{Sch}) and  (\ref{quotient}), 
\begin{align}\label{eight}
(M\quotient{\pi})_{ij}\nonumber
&=\alpha  A(X\quotient{\pi})_{ij}
+(D\quotient{\pi})_{ij}\\
&=
\begin{cases}
   \alpha  k_i+d_i & \text{if }i= j,\\
\alpha  A(H)_{ij}n_j & \text{if }i\neq j. 
\end{cases}
\end{align}
 Hence, letting 
\begin{align}
  \rho_{ij}=
\begin{cases}
   \frac{1}{n_i}\left(x-\alpha  k_i-d_i\right) & \text{if }i= j,\\
 -\alpha  A(H)_{ij} & \text{if }i\neq j, 
\end{cases}
\end{align} 
 we have \quad
$xI-M\quotient{\pi}=(\rho_{ij}n_j).$  
\quad Since
\begin{align*}
(\rho_{ij})&=-\alpha  A(H)+\DD(x), 
\text{ here }\DD(x)_{ij}=\delta_{ij}\frac{1}{n_i}
\left(x-\alpha  k_i-d_i\right),  
\end{align*}
we have (\ref{M}).   \qed \\
\begin{rem}
 We can generalize Corollary \ref{jjoin} to the following Teranishi's result 
 [Ter03, Theorem 3.3]. 

 Let $X=H[X_1,\dots,X_r]$ be a generalized join (composition) 
 and for each $i=1,\dots,r$, let 
$\pi_i=(V_i^1,\dots,V_i^{\alpha_i})$ be an equitable partition of $V(X_i).$ 
Then letting $\pi=\left(\pi_i\right)_{i=1}^r
=\left(\left(V_i^k\right)_{k=1}^{\alpha_i}\right)_{i=1}^r$, 
 $\pi$ is an equitable partition of $V(X)$ and 
\begin{align*}
\charp{A(X)}%
&=\charp{A(X\quotient{\pi})}
\cdot\prod_{i=1}^r\charp{A(X_i\ssetminus \pi_i)}\\
&= \charp{A(X\quotient{\pi})}
\cdot\prod_{i=1}^r\frac{\charp{A(X_i)}}{\charp{A(X_i/\pi_i)}}.
\end{align*}
 When  $X_i$ is $k_i$-regular, letting $\pi_i=\overline{\pi}_i=\{V(X_i)\}$ 
the trivial parition,  $\pi_i$ is an equitable partition of $V(X_i)$ 
 with $\charp{A(X_i/\pi_i)}=x-k_i$, we get 
 (\ref{Xover}) for $\alpha=1, D=0$, and 
similarly for all $\alpha$ and $D=(d_iI_i).$ 
The same proof also applies  in the proof of Corollary \ref{jjoin}, 
 considering the parition $\pi_i$ 
 instead of the trivial partition $\overline{\pi}_i$. 
\end{rem}
\begin{rem}
 As the notation in [CC17, Theorem 3.1], 
\begin{align}\label{zero}
\nonumber \gamma_i& :=1-tk_i+(1-u)(k_i+N_i-1+u)t^2\\ \nonumber
&= -tk_i+1-(1-u)^2t^2+(1-u)t^2(k_i+N_i) 
= -tk_i+s_1+s_2d_i(X)\\
&=-t k_i+\dzi,
\end{align}
 here $N_i=\sum_{j\neq i}A(H)_{ij}n_j.$
 Let $M=-t A(X)+\Dz$. Letting $\alpha=-t,D=D^Z(X)$ in (\ref{eight}) 
 we have
  $(M\quotient{\pi})_{ij}=
 \begin{cases}
\gamma_i     & \text{if }i= j,\\
-t \rho_{ij}n_j & \text{if }i\neq j, 
\end{cases}$

 here $\rho_{ij}=A(H)_{ij}$. 
Let 
$ N_{ij}=
 \begin{cases}
\gamma_i     & \text{if }i= j,\\
-t \sqrt{n_in_j}\rho_{ij} & \text{if }i\neq j. 
\end{cases}
$ 

Then, by letting 
$ \rho'_{ij}=
 \begin{cases}
\frac{\gamma_i}{n_i}     & \text{if }i= j,\\
-t \rho_{ij} & \text{if }i\neq j, 
\end{cases}
$ 

we have 
$M\quotient{\pi}=(\rho'_{ij}n_j)$ and 
$N=(\rho'_{ij}\sqrt{n_in_j}).$ 
 Then 
\begin{align}\label{oneone}
 \det M\quotient{\pi}=
\det \left(-t A(X\quotient{\pi})+(\delta_{ij}d_i^Z(X))\right)
=\det N.
\end{align}
Let \begin{align*}
 M_{X_i}(F)&=I_i-tA(X_i)+(1-u)\left(F-(1-u)I_i\right)t^2\\
&= -t  A(X_i)+s_1I_i+s_2F.
\end{align*} 
 Since 
 $D(X_i)+N_iI_i=(k_i+N_i)I_i=d_i(X)I_i$, we have 
\begin{align}\label{saigo}
 M_{X_i}(D(X_i)+N_iI_i)
&=-t  A(X_i)+\left(s_1+s_2d_i(X)\right)I_i
=-t  A(X_i)+d_i^Z(X)I_i.
\end{align}
 So, by (\ref{gjoin}),(\ref{zero}),(\ref{oneone}),(\ref{saigo}), we have 
\begin{align*}
 Z(u,t)^{-1}
&= (1-(1-u)^2t^2)^{m-n}\det N\cdot \prod_{i=1}^r
 \frac{\det\left(M_{X_i}\left(D(X_i)+N_iI_i\right)\right)}{\gamma_i},
\end{align*}
 which is 
[CC17, Theorem 3.1].
\end{rem}
\bigskip
\hspace{-4mm}{\bf Acknowledgement}\quad 
 I would like to thank Yusuke Ide, Norio Konno, Masashi Kosuda, 
 Hideo Mitsuhashi and Iwao Sato 
 for their help writing this paper.  
 I am indebted to the anonymous referee for  valuable comments 
 and helpful suggestions, 
which led to an improvement of the original manuscript. 
\if01
 The authors are grateful to the unnamed referee for their constructive 
 criticism of the first draft. 
helpful comments and suggestions regarding this paper.
 for his valuable feedback. 
\fi
\vspace*{12mm}
\begin{center}
{\large\bf{References}}
\end{center}
\vspace*{12mm}

\hspace*{-5.8mm}
\bbf{[Bar99] L. Barthordi}  Counting paths in graphs. 
\EM {\bf{45}} (1999) 83-131.
\medskip\\\bbf{[Bas92] H. Bass}
The Ihara-Selberg zeta function of a tree lattice. 
 Internat. J. Math. {\bf{3}} (1992) 717-797.
\medskip\\
\bbf{[BFW16] W. Barrett, A. Francis, and B. Webb}
 Equitable decompositions of graphs with symmetry.
Linear Algebra and its Applications, {\bf 513} 409-434, 2017.
\medskip\\  
\bbf{[CC17] H. Chen and Y. Chen}
Bartholdi zeta functions of generalized join graphs. 
Graphs and Combinatorics {\bf{27}}  (2017).
\if01
\medskip\\
\References{CDR07}{D. M. Cardoso, C. Delorme and P. Rauma}{Laplacian eigenvectors and eigenvalues and almost equitable partitions}{European Journal of Combinatorics}{28}{665-673}{2007}
\medskip\\
\References{CDS79}{D. M. Cvetkovi\'c,M.Doob and H. Sachs}{Spectra of Graphs}
{Academic Press}{}{}{1988} 
\medskip\\
\References{CFM13}{D. M. Cardoso, M. A. A. de Freitas and E. A. Martins}{Spectra of graphs obtained  by a generalization of the join graph operation}{Discrete Mathematics}{313}{733-741}{2013}
\fi
\medskip\\
\References{DSW07}{A. Deng, I. Sato  and Y. Wu}{Homomorphisms, representations 
 and characteristic polynomials of digraphs}
{Linear Algebra and its Applications}{423}{386-407}{2007} 
\medskip\\
\References{DW05}{A. Deng and Y. Wu}{Characteristic polynomials of digraphs 
having a semi-free action}
{Linear Algebra and its Applications}{408}{189-206}{2005} 
\medskip\\
\References{FKL04}{R. Feng, J.H. Kwak and J. Lee}{Characteristic polynomials of graph coverings}{Bulletin of the Australian Mathematical Society}{69}{133-136}{2004}
\medskip\\
\References{FSSW17}{A. Francis, D. Smith, D. Sorensen and B. Webb}{Extensions and applications of equitable decompositions for graph with symmetries}{Linear Algebra and its Applications}{532}{432-462}{2017}
\medskip\\
{[FSW18]} A. Francis, D. Smith and B. Webb, General Equitable Decompositions for Graphs with Symmetries,arXiv:1809.07186 (2018).
\medskip\\\bbf{[FZ99] D. Foata and D. Zeilberger}
A combinatorial proof of Bass's evaluations of the Ihara-Selberg zeta 
function for graphs. 
 Trans. Amer. Math. Soc.  {\bf{351}} (1999) 2257-2274. 
\medskip\\\bbf{[GR01] C. Godsil and G. Royle}  
 Algebraic Graph Theory, Graduate
Texts in Mathematics {\bf 207} Springer, New York, 2001.
\medskip\\\bbf{[GT77] L. Gross and W. Tucker}  
 Generating all graph coverings by permutation voltage assignments.
 Discrete Mathematics, {\bf 18} (1977)  273-283.  
\medskip\\\bbf{[Has89] K., Hashimoto}   
 Zeta functions of finite graphs and representations of $p$-adic 
 groups,  Advanced Studies in Pure Math. {\bf{15}} (1989) 211-280. 
\medskip\\\bbf{[Has90] K., Hashimoto}   
 On the zeta- and $L$-functions of finite graphs,
  International. J.  Math. {\bf{1}} (1990) 381-396. 
\if01
\medskip\\
\References{HOS14}{ M. Hellmuth, L. Ostermeier and P. F. Stadler}  
 {Square property, equitable partitions, and product-like graphs}
{Discrete Mathematics}{320}{92-103}{2014}
\fi
\medskip\\
\References{Iha66}{Y. Ihara}{On discrete subgroups of the two by two projective linear group over {\em p}-Adic fields}{J. Math. Soc. Japan}{18}{219-235}
{1966}
\medskip\\
\References{KL92}{J.H. Kwak and J. Lee}{Characteristic polynomials of graph bundles  $\mathbb{I}$'}
{Linear and Multilinear Algebra}{32}{61-73}{1992} 
\medskip\\
\References{KL01}{J.H. Kwak and Y. S. Kwon}{Characteristic polynomials of graph bundles having voltages in a dihedral group}
{Linear Algebra and its Applications}{336}{99-118}{2001} 
\if01
\medskip\\
\References{KL08}{H. K. Kim and J. Lee}{A generalized characteristic polynomial  of a graph having a semifree action}{Discrete Math.}{308}{555-564}
{2008}
\fi
\medskip\\
\References{KS00}{M. Kotani and T. Sunada}{Zeta functions of finite graphs}
{J. Math. Sci. Univ. Tokyo}{7}{7-25}{2000} 
\medskip\\
\References{MS95}{H. Mizuno and I. Sato}{Characteristic polynomials of some graph covering}
{Discrete Math.}{142}{295-298}{1995} 
\medskip\\
\References{MS97}{H. Mizuno and I. Sato}{Characteristic polynomials of some covers of symmetric digraphs}
{Ars Combinatoria}{45}{3-12}{1997} 
\medskip\\
\References{MS01}{H. Mizuno and I. Sato}{Zeta functions of digraphs}
{Linear Algebra and its Applications}{336}{181-190}{2001} 
\medskip\\
\References{MS02}{H. Mizuno and I. Sato}{Weighted zeta functions of digraphs}
{Linear Algebra and its Applications}{355}{35-48}{2002} 
\medskip\\
\References{MS04}{H. Mizuno and I. Sato}{Weighted zeta functions of graphs}
{Journal of Combinatorial Theory, Series B}{91}{169-183}{2004} 
\medskip\\
\References{Sat99}{I. Sato}{Graph covering and its generalization}
{Yokohama Mathematical Journal}{47}{67-88}{1999} 
\medskip\\
\References{Sat06}{I. Sato}{Weighted Zeta Functions of Graph Coverings}
{The Electronic Journal of Combinatorics}{13}{}{2006} 
\medskip\\
\References{Sat07}{I. Sato}{A New Zeta Function of a Graph}
{International Journal of Algebra}{6}{269-281}{2007} 
\medskip\\
\bbf{[Sch74] A. J. Schwenk} Computing the characteristic polynomial of a graph. In Graphs and Combinatorics, Springer Lecture Notes in Math. 
{\bf{406}} (1974) 153-172.
\if01
\References{Sch74}{A. J. Schwenk}{Computing the characteristic polynomial of a graph}{In Graphs and Combinatorics}{Volume 406 of Lect. Notes Math.}{pages 153-172, Heidelberg, D}{1974}
\fi
\medskip\\\bbf{[Ser80] J. -P. Serre} Trees, Springer-Verlag, New York, 1980. 
\medskip\\\bbf{[ST96] H. M. Stark and A. A. Terras}
Zeta Functions of Finite Graphs and Coverings. 
 Adv. Math.  {\bf{121}} (1996) 124-165. 
\medskip\\
\bbf{[Sun86a] T. Sunada} $L$-functions in geometry and some applications, 
  Springer Lecture Notes in Math. {\bf{1201}} (1986) 
266-284. 
\medskip\\\bbf{[Sun86b] T. Sunada}
Fundamental groups and Laplacian, Proc. of Taniguchi Sympos., 
 Geometry and Analysis on Manifolds,  Springer Lecture Notes in Math. 
{\bf{1339}} (1986) 
248-277.
\medskip\\
\References{Ter03}{Y. Teranishi}{Equitable switching and spectra of graphs}
{Linear Algebra and its Applications}{359}{121-131}{2003} 
\fi
\end{document}